\documentclass{amsart}
\usepackage{amssymb, latexsym}

\def\bi{\begin{itemize}}
\def\bs{\begin{split}}
\def\es{\end{split}}
\def\ba{\begin{align}}
\def\bas{\begin{align*}}
\def\ea{\end{align}}
\def\eas{\end{align*}}
\def\Im{{\hbox{Im}}}

\def\Re{{\hbox{Re}}}

\def\C{{\mathbb C}} % Jim changed to mathbb
\def\R{{{\mathbb R}}}
\def\Z{{{\mathbb Z}}}

\def\O{{{\mathcal O}}}

\def\emph#1{{\it #1}}
\def\textbf#1{{\bf #1}}

\newcommand{\lqlr}{{L^q_t L^r_x}}
\newcommand{\ir}{{I \times \R^4}}
\newcommand{\rr}{{\R \times \R^4}}
\newcommand{\lli}{{\lqlr (\ir)}}
\newcommand{\llr}{{\lqlr (\rr)}}
\newcommand{\sz}{{\dot S^0}}
\newcommand{\szi}{{\sz (\ir)}}

\newcommand{\so}{{\dot S^1}}
\newcommand{\soi}{{\so (\ir)}}
\newcommand{\sor}{{\so (\rr)}}
\newcommand{\st}{{\dot S^2}}

\newcommand{\sk}{{\dot S^k}}
\newcommand{\ski}{{\sk (\ir)}}

\newcommand{\hk}{{\dot H^k}}
\newcommand{\ho}{{\dot H^1}}
\newcommand{\htwo}{{\dot H^2}}
\newcommand{\ls}{{L^6_{t,x}}}
\newcommand{\lsi}{{\ls(\ir)}}
\newcommand{\lsr}{{\ls(\rr)}}
\newcommand{\is}{{I_* \times \R^4}}
\newcommand{\lsis}{{\ls(\is)}}
\newcommand{\sois}{{\so (\is)}}
\newcommand{\szis}{{\sz (\is)}}
\newcommand{\stis}{{\st (\is)}}

\newcommand{\eps}{{\varepsilon}}
\newcommand{\uhi}{{u_{hi}}}
\newcommand{\ulo}{{u_{lo}}}
\newcommand{\uhip}{{u_{hi'}}}
\newcommand{\ecrit}{{E_{crit}}}
\newcommand{\pl}{{P_{lo}}}
\newcommand{\ph}{{P_{hi}}}
\newcommand{\pmed}{{P_{med}}}
\newcommand{\propagate}{{e^{i(t-t_0)\Delta}}}
\newcommand{\propagateo}{{e^{i(t-t_1)\Delta}}}

\theoremstyle{plain}
\newtheorem{theorem}{Theorem}
\newtheorem{proposition}[theorem]{Proposition}
\newtheorem{lemma}[theorem]{Lemma}
\newtheorem{corollary}[theorem]{Corollary}
\theoremstyle{definition}

\newtheorem{definition}[theorem]{Definition}
\newtheorem{remark}[theorem]{Remark}

\numberwithin{equation}{section} \numberwithin{theorem}{section}

\begin{document}

\title[Global well-posedness for 4D energy-critical NLS]%
{Global well-posedness and scattering for the defocusing energy-critical nonlinear Schr\"odinger equation in $\R^{1+4}$}
\author{E. ~Ryckman}
\address{Department of Mathematics\\UCLA\\Los Angeles, CA 90095}
\author{M.~Visan}
\address{Department of Mathematics\\UCLA\\Los Angeles, CA 90095}
%\subjclass{35Q55} \keywords{nonlinear Schr\"odinger equation,
%well-posedness}

\vspace{-0.3in}
\begin{abstract}
We obtain global well-posedness, scattering, uniform regularity, and global $L^6_{t,x}$ spacetime bounds for energy-space
solutions to the defocusing energy-critical nonlinear Schr\"odinger equation in $\R\times\R^4$.  Our arguments closely follow those in
\cite{ckstt:gwp}, though our derivation of the frequency-localized interaction Morawetz estimate is somewhat simpler.
As a consequence, our method yields a better bound on the $L^6_{t,x}$-norm.
\end{abstract}

\maketitle

\section{Introduction}

We study the following initial value problem for the cubic defocusing nonlinear Schr\"odinger equation in $\R\times\R^4$
\begin{equation}\label{schrodinger equation}
\begin{cases}
i u_t +\Delta u = |u|^2 u\\
u(0,x) = u_0(x)
\end{cases}
\end{equation}
where $u(t,x)$ is a complex-valued function in spacetime $\R_t\times \R^4_x$.

This equation has the Hamiltonian
\begin{equation}\label{energy}
E(u(t))=\int \frac{1}{2}|\nabla u (t,x)|^2 + \frac{1}{4}|u(t,x)|^4dx.
\end{equation}
Since \eqref{energy} is preserved by the flow corresponding to \eqref{schrodinger equation} we shall refer to it
as the \emph{energy} and often write $E(u)$ for $E(u(t))$.

A second conserved quantity we will occasionally rely on is the mass $\|u(t)\|_{L^2_x(\R^4)}^2$.  However, since
the equation is $L_x^2$-supercritical with respect to the scaling (see below), we do not have bounds on the mass that
are uniform across frequencies (indeed, the low frequencies may simultaneously have small energy and large mass).

We are primarily interested in the cubic defocusing equation \eqref{schrodinger equation} since it is critical
with respect to the energy norm.  That is, the scaling $u \mapsto u^\lambda$ where
\begin{equation}\label{scaling}
u^\lambda (t,x) :=
\frac{1}{\lambda}u(\frac{t}{\lambda^2},\frac{x}{\lambda}),
\end{equation}
maps a solution to \eqref{schrodinger equation} to another solution to \eqref{schrodinger equation}, and $u$ and
$u^\lambda$ have the same energy.

It is known that if the initial data $u_0$ has finite energy, then
\eqref{schrodinger equation} is locally well-posed (see, for
instance \cite{cwI}). That is, there exists a unique local-in-time
solution that lies in $C_t^0 \ho_x \cap L^6_{t,x}$, and the map from
the initial data to the solution is locally Lipschitz in these
norms. If in addition the energy is small, it is known
that the solution exists globally in time and scattering occurs;
that is, there exist solutions $u_\pm$ of the free Schr\"odinger
equation $(i \partial_t +\Delta)u_\pm = 0$ such that $\|u(t) - u_\pm
(t)\|_{\ho_x} \rightarrow 0$ as $t \rightarrow \pm \infty$. However,
for initial data with large energy, the local well-posedness
arguments do not extend to give global well-posedness.

Global well-posedness in $\dot{H}^1_x(\R^3)$ for the energy-critical NLS in the case of
large finite-energy, radially-symmetric initial data was first obtained by Bourgain (\cite{borg:scatter}, \cite{borg:book})
and subsequently by Grillakis \cite{grillakis:scatter}. Tao \cite{tao: gwp radial} settled the problem for arbitrary
dimensions (with an improvement in the final bound due to a simplification of the argument), but again only for radially
symmetric data. A major breakthrough in the field was made by J. Colliander, M. Keel, G. Staffilani, H. Takaoka, and T. Tao in
\cite{ckstt:gwp} where they obtained global
well-posedness for the energy-critical NLS in dimension $n=3$ with arbitrary data. In dimensions $n \geq 4$, the problem was
still open.

The main result of this paper is the global well-posedness statement for
\eqref{schrodinger equation} in the energy space,

\begin{theorem}\label{lemma THE THEOREM}
For any $u_0$ with finite energy $E(u_0) < \infty$, there exists a
unique global solution $u \in C_t^0 \ho_x \cap L^6_{t,x}$ to
\eqref{schrodinger equation} such that
\begin{equation}\label{l6 bounds}
\int_{-\infty}^{\infty} \int_{\R^4} |u(t,x)|^6 dx dt \leq C(E(u_0))
\end{equation}
for some constant $C(E(u_0))$ depending only on the energy.
\end{theorem}

\subsection{Outline of the proof of Theorem \ref{lemma THE THEOREM}}
We first notice that it suffices to prove
Theorem \ref{lemma THE THEOREM} for Schwartz functions. Indeed, if one obtains a uniform $L^6_{t,x}(I\times\R^4)$
bound for all Schwartz solutions and all compact intervals $I$, one can approximate arbitrary finite-energy initial
data by Schwartz initial data and use perturbation theory to prove that the corresponding sequence of solutions to
\eqref{schrodinger equation} converges in $\dot{S}^1(I\times\R^4)$ to a finite-energy solution of
\eqref{schrodinger equation}.  See Sections 1.2 and 2.1 for the notation
and definitions appearing in the outline of the proof.

For an energy $E \geq 0$ we define the quantity $M(E)$ by
$$M(E):=\sup\|u\|_{L^6_{t,x}(I\times\R^4)},$$
where $I \subset \R$ ranges over all compact time intervals, and
$u$ ranges over all Schwartz solutions to \eqref{schrodinger
equation} on $I\times \R^4$ with $E(u)\leq E$.  For $E<0$ we define
$M(E)=0$ since, of course, there are no negative energy solutions.  Our task is to show that $M(E) < \infty$ for all $E$.

Let us note that by the small-energy well-posedness result discussed above, we know that $M(E)$ is finite for $E$
sufficiently small.
We will assume for a contradiction that $M(E)$ is not always finite.  {}From perturbation theory (see Lemma \ref{lemma long time})
it follows that the set $\{E : M(E)<\infty\}$ is open.  Since it is also connected and contains zero, there exists a
\emph{critical energy} $0 < E_{crit} < \infty$ such that $M(E_{crit}) = \infty$ but
$M(E) < \infty$ for all $E < E_{crit}$.  From the definition of $E_{crit}$ and the $\ls$
well-posedness theory (see Section 3 for details) we get:

\begin{lemma}[Induction on energy hypothesis]\label{lemma induct on energy}
Let $t_0 \in \R$ and let $v(t_0)$ be a Schwartz function with
$E(v_0) \leq E_{crit}-\eta$ for some $\eta > 0$.  Then there
exists a global Schwartz solution $v$ of \eqref{schrodinger
equation} on $\rr$ with initial data $v(t_0)$ at time $t_0$, such
that
$$\|v\|_{\ls(\rr)} \leq M(E_{crit} - \eta) .$$
Moreover we have
$$\|v\|_{\sor} \leq C(E_{crit}, M(E_{crit} - \eta)).$$
\end{lemma}

We will need a few small parameters for the contradiction argument
(albeit two less than were necessary in \cite{ckstt:gwp}).
Specifically we will need
$$1 \gg \eta_0 \gg \eta_1 \gg \eta_2 \gg \eta_3 \gg \eta_4 > 0$$ where each
$\eta_j$ is allowed to depend on the critical energy and any of
the larger $\eta$'s. We will choose $\eta_j$ small enough that, in
particular, it will be smaller than any constant depending on the
previous $\eta$'s used in the argument.

As $M(E_{crit})$ is infinite, given any $\eta_4> 0$ there exist a compact interval $I_* \subset \R$ and a Schwartz
solution $u$ to \eqref{schrodinger equation} on $\is$ with
$ E(u) \leq E_{crit}$ but
\begin{equation}\label{l6 HUGE}
\|u\|_{\ls (\is)} > 1/\eta_4.
\end{equation}
Note that we may assume $E(u)\geq \frac{1}{2}E_{crit}$, since otherwise we would get
$$\|u\|_{\ls (\is)} \leq M(\frac{1}{2}E_{crit})<\infty$$ and we would be done.

This suggests we make the following definition.

\begin{definition}
A \emph{minimal energy blowup solution} of \eqref{schrodinger
equation} is a Schwartz solution $u$ on a time interval $I_*
\subset \R$ with energy
\begin{equation}\label{minimal energy blowup solution}
\frac{1}{2}E_{crit} \leq E(u(t)) \leq E_{crit}
\end{equation}
and huge $\ls$-norm in the sense of \eqref{l6 HUGE}.
\end{definition}

%We will derive the desired contradiction by showing that $\|u\|_{\ls (\is)}\leq C(\eta_0,\eta_1,\eta_2,\eta_3)$.

Note that \eqref{minimal energy blowup solution} implies the kinetic energy bound
\begin{equation}\label{kinetic energy bound}
\|u\|_{L^\infty_t \ho_x (\is)} \sim 1,
\end{equation}
while from Sobolev embedding we obtain a bound on the potential energy,
\begin{equation}\label{potential energy bound}
\|u\|_{L^\infty_t L^4_x (\is)} \lesssim 1.
\end{equation}

In Sections 2 and 3 we recall the Strichartz estimates and the
perturbation theory we will use throughout the proof of Theorem
\ref{lemma THE THEOREM}.  Many of the ideas of these sections have been previously developed,
 and in a few cases we content ourselves with citing the
relevant source (e.g., Lemmas \ref{lemma linear strichartz} and
\ref{lemma bilinear strichartz}).

We expect that a minimal energy blowup solution should be localized in both physical and frequency space.  For if not,
it could be decomposed into two essentially separate solutions, each with strictly smaller energy than the original.
By Lemma~\ref{lemma induct on energy} we can then extend these smaller energy solutions to all of $I_*$.  As each of the
separate evolutions exactly solves \eqref{schrodinger equation}, we expect their sum to solve \eqref{schrodinger equation}
approximately.  We could then use the perturbation theory results and the bounds from Lemma
\ref{lemma induct on energy} to bound $\|u\|_{L^6_{t,x}}$ in terms of $\eta_0, \eta_1, \eta_2,\text{ and } \eta_3$,
thus contradicting the fact that $\eta_4$ can be chosen arbitrarily small.

This intuition will underpin the frequency localization argument
we give in Section 4.  The spatial concentration result follows in
a similar manner, but is a bit more technical. For instance,
we restrict our analysis to a subinterval
$I_0 \subset I_*$ and will need to use both frequency localization
and the fact that the potential energy of a minimal energy blowup solution
is bounded away from zero.

In Section 5 we obtain the frequency-localized Morawetz inequality \eqref{flim}, which will be used
to derive a contradiction to the frequency localization and spatial concentration results just described.

A typical example of a Morawetz inequality for \eqref{schrodinger equation} (see \cite{linstrauss}) is the bound
$$
\int_I \int_{\R^4}\frac{|u(t,x)|^4}{|x|}dxdt\lesssim \sup_{t\in I}\|u(t)\|_{\dot{H}^{1/2}(\R^4)}^2,
$$
for all time intervals $I$ and all Schwartz solutions $u:I\times\R^4\to \C$.

This estimate is not particularly useful for the energy-critical problem since the $\dot{H}^{1/2}_x$ norm is
supercritical with respect to the scaling \eqref{scaling}. To get around this problem, Bourgain and Grillakis introduced
a spatial cutoff obtaining the variant
$$
\int_I \int_{|x|\leq A|I|^{1/2}}\frac{|u(t,x)|^4}{|x|}dxdt\lesssim A|I|^{1/2} E(u),
$$
for all $A\geq 1$, where $|I|$ denotes the length of the time interval $I$. While this estimate is better suited for
the critical NLS (it involves the energy on the right-hand side), it only prevents concentration of $u$ at the spatial
origin $x=0$. This is especially useful in the spherically-symmetric case $u(t,x)=u(t,|x|)$, since the spherical symmetry
combined with the bounded energy assumption can be used to show that $u$ cannot concentrate at any location but the
spatial origin. However, it does not provide much information about the solution away from the origin.
Following \cite{ckstt:gwp}, we circumvent this problem by using a frequency-localized interaction Morawetz inequality.

While the previously mentioned Morawetz inequalities were \textit{a priori} estimates, the frequency-localized
interaction Morawetz inequality we will develop is not; it only applies to minimal
energy blowup solutions. Our model in obtaining the frequency-localized interaction Morawetz estimate is \cite{ckstt:gwp}.
However, our argument is not as technical as theirs, lacking the need for spatial localization. It is this
simplification that will yield an improvement in the final bound on the $L^6_{t,x}$-norm.

A corollary of \eqref{flim} is good $L^3_{t,x}$ control over the high-frequency part of a minimal energy blowup solution.
One then has to use Sobolev embedding to bootstrap this $L^3_{t,x}$ control to $L^6_{t,x}$ control. However, one needs
to make sure that the solution is not shifting its energy from low to high frequencies causing the $L^6_{t,x}$-norm
to blow up while the $L^3_{t,x}$-norm stays bounded. This is done in Section~6, where we prove a frequency-localized
mass estimate that prevents energy evacuation to high modes.  We put all these pieces together in Section~7
where the contradiction argument is concluded.

In Section 8 we comment on the tower bound we get on the $L^6_{t,x}$-norm in Theorem \ref{lemma THE THEOREM}, and we
show how this bound yields scattering, asymptotic completeness, and uniform regularity.

As will certainly become clear to the reader, our paper relies heavily on the arguments developed in
\cite{ckstt:gwp}. One should note a few differences though, mainly related to the Strichartz norms and the
frequency-localized interaction Morawetz inequality. While it is true that in higher dimensions one has more
Strichartz estimates, we lack
control over $L_t^2L_x^{\infty}$ (for which one gets a logarithmic divergence). This turns out not to be a problem most
of the time, since the triangles in which the low- and high-frequency parts of the minimal energy blowup solution live are
large enough in four dimensions. (Indeed, by Bernstein, the low-frequency portion of the solution has finite spacetime norm for every
$L_t^p L_x^q$ in the closed triangle with vertices $L_t^2 L_x^\infty$, $L_t^\infty L_x^4$, and $L_{t,x}^\infty$,
except for the vertex $L_t^2 L_x^\infty$.  By interpolation, the high-frequency portion has finite spacetime norm for every $L_t^p L_x^q$ in the
closed triangle with vertices $L_{t,x}^3$, $L_t^\infty L_x^2$, and $L_t^\infty L_x^4$.)  However, when it comes to controlling the error terms generated
by the frequency localization in the interaction Morawetz inequality, we have to do something different. While in \cite{ckstt:gwp} the
authors strive to gain better control on low frequencies, we will use multilinear operator theory
(specifically a theorem of Coifman and Meyer \cite{coifmey:1}, \cite{coifmey:2}) to gain control
over the high frequencies. Being able to control the high frequencies in $L^3_{t,x}$ will save us from having to localize
the interaction Morawetz inequality in space as well. As a consequence, our argument is somewhat simpler and it yields a smaller
tower bound on the $L^6_{t,x}$-norm.

\textbf{Acknowledgments}:
We thank Terence Tao for valuable comments on this paper and explanatory details related to \cite{ckstt:gwp}.

\subsection{Notation}
We will often use the notation $X \lesssim Y$ whenever there exists some constant $C$, possibly depending on the
critical energy but not on any other parameters, so that $X \leq CY$.
Similarly we will use $X \sim Y$ if $X \lesssim Y \lesssim X$.  We use $X \ll Y$ if $X \leq cY$ for some small constant
$c$, again possibly depending on the critical energy. We will use the notation $X+ := X+\epsilon$, for some
$0<\epsilon \ll 1$; similarly $X- := X-\epsilon$. We also use the notation $\langle x\rangle :=(1+|x|^2)^{1/2}.$

We define the Fourier transform on $\R^4$ to be $$\hat f(\xi) :=
\int_{\R^4} e^{-2 \pi i x \cdot \xi} f(x) dx.$$

We will make frequent use of the fractional differentiation
operators $|\nabla|^s$ defined by $$\widehat{|\nabla|^sf}(\xi) :=
|\xi|^s \hat f (\xi).$$  These define the homogeneous Sobolev
norms
$$\|f\|_{\dot H^s_x} := \| |\nabla|^s f \|_{L^2_x (\R^4)}.$$

Let $e^{it\Delta}$ be the free Schr\"odinger propagator.  In physical space this is given by the formula
$$e^{it\Delta}f(x) = \frac{1}{(4 \pi i t)^2} \int_{\R^4} e^{i|x-y|^2/4t} f(y) dy,$$
while in frequency space one can write this as
\begin{equation}\label{fourier rep}
\widehat{e^{it\Delta}f}(\xi) = e^{-4 \pi^2 i t |\xi|^2}\hat
f(\xi).
\end{equation}

In particular, the propagator preserves the above Sobolev norms
and obeys the \emph{dispersive inequality}
\begin{equation}\label{dispersive ineq}
\|e^{it\Delta}f\|_{L^\infty_x(\R^4)} \lesssim
|t|^{-2}\|f\|_{L^1_x(\R^4)}
\end{equation}
for all times $t$.  We also recall \emph{Duhamel's formula}
\begin{align}\label{duhamel}
u(t) = e^{i(t-t_0)\Delta}u(t_0) - i \int_{t_0}^t e^{i(t-s)\Delta}(iu_t + \Delta u)(s) ds.
\end{align}

We will use the notation $\O(X)$ to denote a quantity that resembles $X$; that is a finite linear combination of
terms that look like $X$, but possibly with some factors replaced by
their complex conjugates.  For example we will write
\begin{equation}\label{slash o notation}
|u+v|^2(u+v) = \sum_{j=0}^3 \O(u^j v^{3-j}).
\end{equation}

We will occasionally use subscripts to denote spatial derivatives
and will use the summation convention over repeated indices.

We will also need some Littlewood-Paley theory.  Specifically, let
$\varphi(\xi)$ be a smooth bump adapted to the ball $|\xi| \leq 2$
equalling one on the ball $|\xi| \leq 1$.  For each dyadic number
$N \in 2^\Z$ we define the Littlewood-Paley operators
\begin{align*}
\widehat{P_{\leq N}f}(\xi) &:=  \varphi(\xi/N)\hat f (\xi)\\
\widehat{P_{> N}f}(\xi) &:=  (1-\varphi(\xi/N))\hat f (\xi)\\
\widehat{P_N f}(\xi) &:=  (\varphi(\xi/N) - \varphi (2 \xi /N))
\hat f (\xi).
\end{align*}
Similarly we can define $P_{<N}$, $P_{\geq N}$, and $P_{M < \cdot
\leq N} := P_{\leq N} - P_{\leq M}$, whenever $M$ and $N$ are
dyadic numbers.  We will frequently write $f_{\leq N}$ for
$P_{\leq N} f$ and similarly for the other operators.

The Littlewood-Paley operators commute with derivative operators,
the free propagator, and the conjugation operator.  They are self-adjoint
and bounded on every $L^p$ and $\dot H^s$ space for $1 \leq p \leq
\infty$ and $s\geq 0$.  They also obey the following Sobolev and Bernstein
estimates
\begin{align*}
\|P_{\geq N} f\|_{L^p} &\lesssim N^{-s} \| |\nabla|^s P_{\geq N} f \|_{L^p},\\
\| |\nabla|^s  P_{\leq N} f\|_{L^p} &\lesssim N^{s} \| P_{\leq N} f\|_{L^p},\\
\| |\nabla|^{\pm s} P_N f\|_{L^p} &\sim N^{\pm s} \| P_N f \|_{L^p},\\
\|P_{\leq N} f\|_{L^q} &\lesssim N^{\frac{4}{p}-\frac{4}{q}} \|P_{\leq N} f\|_{L^p},\\
\|P_N f\|_{L^q} &\lesssim N^{\frac{4}{p}-\frac{4}{q}} \| P_N f\|_{L^p},
\end{align*}
whenever $s \geq 0$ and $1 \leq p \leq q \leq \infty$.

For instance, we can use the above Bernstein estimates and
\eqref{kinetic energy bound} to bound the mass of high frequencies
\begin{equation}\label{mass high freq bound}
\|P_{>M}u\|_{L^2(\R^4)} \lesssim \frac{1}{M} \text{ for all }M \in
2^{\Z}.
\end{equation}

%%%%%%%%%%%%%%%%%%%%%%%%%%%%%%%%%%%%%%%%%%%%%%%%%%%%%%%%%%%%%%%%%%%%%%%%%%%%%%%%%%%%%%%%%%%
%
%
%                                   Section
%
%
%%%%%%%%%%%%%%%%%%%%%%%%%%%%%%%%%%%%%%%%%%%%%%%%%%%%%%%%%%%%%%%%%%%%%%%%%%%%%%%%%%%%%%%%%%%

\section{Strichartz numerology}
In this section we recall Strichartz and bilinear Strichartz estimates in $\R\times\R^4$ and develop trilinear Strichartz estimates.

We use $\lqlr$ to denote the spacetime norm
$$
\|u\|_{\llr} :=\Bigl(\int_{\R}\Bigl(\int_{\R^4} |u(t,x)|^r dx \Bigr)^{q/r} dt \Bigr)^{1/q},
$$
with the usual modifications when $q$ or $r$ is infinity, or when the domain $\R \times \R^4$ is replaced by some
smaller spacetime region.  When $q=r$ we abbreviate $\lqlr$ by $L^q_{t,x}$.

\subsection{Linear and bilinear Strichartz estimates}\label{strichartz section}
We say that a pair of exponents $(q,r)$ is \emph{admissible} if $\tfrac{2}{q} + \tfrac{4}{r} = 2$ and
$2 \leq q,r \leq \infty$. If $I \times \R^4$ is a spacetime slab, we define the $\szi$ \emph{Strichartz norm} by
\begin{equation}\label{s0}
\|u\|_{\szi} := \sup (\sum_N \| P_N u \|^2_{\lli})^{1/2}
\end{equation}
where the $\sup$ is taken over all admissible pairs $(q,r)$.  For
$k = 1,2$ we define the $\dot
S^k (\ir)$ \emph{Strichartz norm} by
$$\|u\|_{\dot S^k (\ir)} := \| \nabla ^k u \|_{\szi}.$$

We observe the inequality
\begin{equation}\label{square sum}
\|(\sum_N |f_N|^2 )^{1/2}\|_{\lli} \leq (\sum_N
\|f_N\|^2_{\lli})^{1/2}
\end{equation}
for all $2 \leq q,r \leq \infty$ and arbitrary functions $f_N$, which one proves by interpolating between the trivial
cases $(2,2)$, $(2,\infty)$, $(\infty,2)$, and $(\infty,\infty)$. In particular, \eqref{square sum} holds for all
admissible exponents $(q,r)$.  Combining this with the Littlewood-Paley inequality we find
\begin{align*}
\| u \|_{\lli}& \lesssim \|(\sum_N |P_N u|^2)^{1/2}\|_{\lli}\\
            & \lesssim (\sum_N \|P_N u \|^2_{\lli})^{1/2}\\
            & \lesssim \| u \|_{\szi},
\end{align*}
which in particular implies
\begin{equation}\label{grad less s1}
\|\nabla u \| _{\lli} \lesssim \|u\|_{\soi}.
\end{equation}

In fact, by \eqref{grad less s1} and Sobolev embedding the $\so$
norm controls the following spacetime norms:

\begin{lemma}\label{lemma strichartz norms}
For any Schwartz function $u$ on $\ir$, we have
\begin{multline}\label{strichartz norms}
\|\nabla u\|_{L^\infty_t L^2_x} + \|\nabla u \|_{L^6_t L^{12/5}_x}
+ \|\nabla u\|_{L^4_t L^{8/3}_x} + \|\nabla u\|_{L^3_{t,x}} +
\|\nabla u\|_{L^2_t L^4_x}\\ + \|u\|_{L^\infty_t L^4_x} +
\|u\|_{L^6_{t,x}} + \|u\|_{L^4_t L^8_x} + \|u\|_{L^3_t L^{12}_x}
\lesssim \|u\|_{\so}
\end{multline}
where all spacetime norms are on $\ir$.
\end{lemma}

One has the following standard Strichartz estimates (see for instance \cite{tao:keel}):

\begin{lemma}\label{lemma linear strichartz}
Let $I$ be a compact time interval, and let $u : \ir \rightarrow
\C$ be a Schwartz solution to the forced Schr\"odinger equation
\begin{equation*}
i u_t + \Delta u = \sum_{m=1}^M F_m
\end{equation*}
for some Schwartz functions $F_1 ,\dots,F_M$.  Then we have
\begin{equation}
\|u\|_{\ski} \lesssim \|u(t_0)\|_{\hk (\R^4)} + \sum_{m=1}^M
\|\nabla^k F_m \|_{L^{q'_m}_t L^{r'_m}_x (\ir)}
\end{equation}
for any integer $k \geq 0$, any time $t_0 \in I$, and any
admissible exponents $(q_1,r_1),\dots,$ $(q_m,r_m)$, where as usual
$p'$ denotes the dual exponent to $p$, that is $1/p + 1/p' = 1$.
\end{lemma}

We also recall without proof the bilinear Strichartz estimates which were obtained in \cite{ckstt:gwp} (see their Lemma 3.4):
\begin{lemma}\label{lemma bilinear strichartz}
Let $n \geq 2$.  For any spacetime slab $I \times \R^n$, any $t_0
\in I$, and any $\delta > 0$, we have
\begin{multline}
\|uv\|_{L^2_{t,x} (I \times \R^n)} \leq
C(\delta)(\|u(t_0)\|_{\dot H_x^{-1/2 +\delta}} + \|(i \partial_t + \Delta)u\|_{L^1_t \dot H_x^{-1/2 + \delta}(I \times \R^n)})\\
\times (\|v(t_0)\|_{\dot H_x^{\frac{n-1}{2}-\delta}} + \|(i
\partial_t + \Delta)v\|_{L^1_t \dot H_x^{\frac{n-1}{2}-\delta}(I \times \R^n)}).
\end{multline}
\end{lemma}

\subsection{Trilinear Strichartz estimates}
We will also need the following estimate:

\begin{lemma}\label{lemma trilinear strichartz}
For $k = 0,1,2$, any slab $\ir$, and any smooth functions
$v_1,v_2,v_3$ on this slab, we have
\begin{equation}\label{trilinear strichartz}
\|\nabla^k \O(v_1 v_2 v_3)\|_{L^1_t L^2_x} \lesssim
\sum_{\{a,b,c\}=\{1,2,3\}} \|v_a\|_{\sk} \min(\|v_b\|_{\so}
\|v_c\|_{L^6_{t,x}}, \|v_b\|_{\ls}
\|v_c\|_{\so}),
\end{equation}
where all spacetime norms are on $\ir$.  Similarly we
have
\begin{equation}\label{trilinear 4/3}
\|\nabla \O (v_1 v_2 v_3)\|_{L^2_t L^{4/3}_x}
\lesssim \sum_{\{a,b,c\}=\{1,2,3\}}\|\nabla v_a\|_{L_t^6L_x^{12/5}} \|v_b\|_{L^6_{t,x}} \|v_c\|_{L_{t,x}^6}
\lesssim \prod_{j=1}^3 \|v_j\|_{\so}.
\end{equation}
\end{lemma}

\begin{proof}
First consider the $k=0$ case of \eqref{trilinear strichartz}. The claim follows from \eqref{strichartz norms} estimating
$$
\|\O ( v_1 v_2 v_3)\|_{L^1_t L^2_x}
  \lesssim \| v_1\|_{L^2_t L^4_x} \|v_2\|_{L^3_t L^{12}_x} \|v_3\|_{L^6_{t,x}}.
$$
Applying the Leibnitz rule we see that the case $k=1$ reduces to obtaining estimates of the form
$$
\|\O ((\nabla v_1) v_2 v_3)\|_{L^1_t L^2_x}
  \lesssim \|\nabla v_1\|_{L^2_t L^4_x} \|v_2\|_{L^3_t L^{12}_x} \|v_3\|_{L^6_{t,x}}.
$$
The lemma again follows from \eqref{strichartz norms}.  The $k = 2$ case of \eqref{trilinear strichartz} proceeds similarly
using estimates such as
$$\|\O ((\nabla^2 v_1) v_2 v_3)\|_{L^1_t L^2_x}
    \lesssim \|\nabla^2 v_1\|_{L^2_t L^4_x} \|v_2\|_{L^3_t L^{12}_x} \|v_3\|_{L^6_{t,x}}
$$
and
$$\|\O ((\nabla v_1) (\nabla v_2) v_3)\|_{L^1_t L^2_x}
    \lesssim \|\nabla v_1\|_{L^2_t L^4_x} \|\nabla v_2\|_{L^3_t L^{12}_x} \|v_3\|_{L^6_{t,x}}.
$$

Finally, the estimate \eqref{trilinear 4/3} follows from \eqref{strichartz norms} and H\"older's inequality,
$$
\|\O ((\nabla v_1) v_2 v_3)\|_{L^2_t L^{4/3}_x}
    \lesssim \|\nabla v_1\|_{L_t^6L_x^{12/5}} \|v_2\|_{L^6_{t,x}} \|v_3\|_{L_{t,x}^6}.
$$
\end{proof}

The following is a variant of the above lemma adapted to the case where some factors are high frequency and others
are low frequency.

\begin{lemma}\label{lemma trilinear strichartz hilo}
Suppose $v_{hi}$ and $v_{lo}$ are functions on $\ir$ such that
\begin{align*}
\|v_{hi}\|_{\sz} + \|(i \partial_t + \Delta)v_{hi}\|_{L^1_t L^2_x}
& \lesssim \varepsilon K\\
\|v_{hi}\|_{\so} + \|\nabla (i \partial_t + \Delta)v_{hi}\|_{L^1_t L^2_x}
& \lesssim  K\\
\|v_{lo}\|_{\so} + \|\nabla (i \partial_t + \Delta)v_{lo}\|_{L^1_t L^2_x}
& \lesssim K\\
\|v_{lo}\|_{\st} + \|\nabla^2 (i \partial_t + \Delta)v_{lo}\|_{L^1_t L^2_x}
& \lesssim \varepsilon K
\end{align*}
for some constants $K>0$ and $0 < \varepsilon \ll 1$ (where all
spacetime norms are on $\ir$).  Then for $j = 1,2$ and any $0 <
\delta \ll 1$, we have
\begin{equation*}
\|\nabla \O (v_{hi}^j v_{lo}^{3-j})\|_{L^2_t L^{4/3}_x (\ir)}
\lesssim \varepsilon^{1-2\delta}K^3.
\end{equation*}
\end{lemma}

\begin{remark}
The importance of this lemma lies in the gain of $\eps^{1-2 \delta}$.  The
$\sz$ bound effectively restricts $v_{hi}$ to high frequencies;
similarly the $\st$ bound restricts $v_{lo}$ to low frequencies.
Comparing the conclusion with \eqref{slash o notation} we see that
the components of the nonlinearity in \eqref{schrodinger equation}
arising from the interaction between high and low frequencies are
rather weak.  While this will be especially important for the
frequency localization result (Proposition \ref{lemma freq loc}),
the idea of controlling frequency interactions will appear
repeatedly.
\end{remark}

\begin{proof}
Throughout, all spacetime norms will be in $\ir$.  We begin by
normalizing $K = 1$.  By Leibnitz we have
\begin{equation*}
\|\nabla \O( v_{hi}^j v_{lo}^{3-j} ) \|_{L^2_t L^{4/3}_x}
   \lesssim \|\O(v_{hi}^j v_{lo}^{2-j} \nabla v_{lo})\|_{L^2_tL^{4/3}_x} + \|\O(v_{hi}^{j-1} v_{lo}^{3-j} \nabla v_{hi})\|_{L^2_t L^{4/3}_x}.
\end{equation*}
Consider the first term.  By H\"older we can bound this by
\begin{align*}
\|\nabla v_{lo}\|_{L^\infty_t L^4_x} \|v_{hi}\|_{L^2_t L^4_x} \|v_{hi}\|_{L^\infty_t L^4_x}^{j-1} \|v_{lo}\|_{L^\infty_t L^4_x}^{2-j}
    &\lesssim \|v_{lo}\|_{\st} \|v_{hi}\|_{\sz} \|v_{hi}\|_{\so}^{j-1} \|v_{lo}\|_{\so}^{2-j}\\
    &\lesssim \varepsilon^2.
\end{align*}

Now consider the second term for the case $j = 2$.  We can bound
\begin{align*}
\|\O(v_{hi} v_{lo} \nabla v_{hi})\|_{L^2_t L^{4/3}_x}
   &\lesssim \|\nabla v_{hi} \|_{L^2_t L^4_x} \|v_{lo}\|_{L^\infty_{t,x}} \|v_{hi}\|_{L^\infty_t L^2_x}\\
   &\lesssim \|v_{hi} \|_{\so} \|v_{lo}\|_{L^\infty_{t,x}} \|v_{hi}\|_{\sz}.
\end{align*}
Decomposing $v_{lo} = \sum_N P_N v_{lo}$ and using Bernstein we get
$$
\|v_{lo}\|_{L^\infty_{t,x}}
    \lesssim (\sum_{N}\|P_N v_{lo}\|_{L^\infty_{t,x}}^2)^{1/2}
    \lesssim (\sum_{N}\|\nabla P_N v_{lo}\|_{L^\infty_tL^4_x}^2)^{1/2}
    \lesssim \|v_{lo}\|_{\st}.
$$
Hence, we obtain the bound
\begin{equation*}
\|\O(v_{hi} v_{lo} \nabla v_{hi})\|_{L^2_t L^{4/3}_x}
    \lesssim \| v_{hi} \|_{\so} \|v_{lo}\|_{\st} \|v_{hi}\|_{\sz}
    \lesssim \varepsilon^2
\end{equation*}
which is acceptable.

Finally consider the second term for the case $j = 1$.  We split the $v_{lo}$ terms dyadically and use H\"older to bound
\begin{align*}
\|\O(v_{lo}^2 \nabla v_{hi})\|_{L^2_t L^{4/3}_x}
   & \lesssim \sum_{N_1 \geq N_2} \|\O((P_{N_1}v_{lo})(P_{N_2}v_{lo})\nabla v_{hi})\|_{L^2_t L^{4/3}_x}\\
   & \lesssim \sum_{N_1 \geq N_2} \|\O((P_{N_2}v_{lo})\nabla v_{hi})\|_{L^2_{t,x}} \|P_{N_1}v_{lo}\|_{L^\infty_t L^4_x}.
\end{align*}
Now $\|P_{N_1}v_{lo}\|_{L^\infty_t L^4_x} \lesssim \|v_{lo}\|_{\so}
\lesssim 1$, while by Bernstein,
\begin{equation*}
\|P_{N_1}v_{lo}\|_{L^\infty_t L^4_x}
\lesssim N_1^{-1} \|P_{N_1} \nabla v_{lo}\|_{L^\infty_t L^4_x}
\lesssim N_1^{-1}\|v_{lo}\|_{\st}
\lesssim \varepsilon N_1^{-1}.
\end{equation*}
So $\|P_{N_1}v_{lo}\|_{L^\infty_t L^4_x} \lesssim \min(1,\varepsilon N_1^{-1})$.

By Lemma \ref{lemma bilinear strichartz} we have
\begin{multline*}
\|\O((P_{N_2}v_{lo})\nabla v_{hi})\|_{L^2_{t,x}} \lesssim
(\|\nabla v_{hi}(t_0)\|_{\dot H_x^{-1/2 + \delta}} + \|(i
\partial_t + \Delta)\nabla v_{hi}\|_{L^1_t \dot
H_x^{-1/2+\delta}})\\
\times (\|P_{N_2}v_{lo}(t_0)\|_{\dot H_x^{3/2-\delta}} + \|(i
\partial_t + \Delta) P_{N_2} v_{lo}\|_{L^1_t \dot
H_x^{3/2-\delta}}).
\end{multline*}
By interpolation we can bound
\begin{equation*}
\|\nabla v_{hi}(t_0)\|_{\dot H_x^{-1/2 + \delta}} + \|(i
\partial_t + \Delta)\nabla v_{hi}\|_{L^1_t \dot
H_x^{-1/2+\delta}} \lesssim \varepsilon^{1/2-\delta},
\end{equation*}
while by Bernstein and the $\so$ and $\st$ bounds we have
\begin{equation*}
\|P_{N_2}v_{lo}(t_0)\|_{\dot H_x^{3/2-\delta}} + \|(i \partial_t + \Delta) P_{N_2} v_{lo}\|_{L^1_t \dot H_x^{3/2-\delta}}
\lesssim \min(N_2^{1/2-\delta}, \eps N_2^{-1/2-\delta}).
\end{equation*}
Putting these together we see that
\begin{equation*}
\|\O(v_{lo}^2 \nabla v_{hi})\|_{L^2_t L^{4/3}_x}
  \lesssim \sum_{N_1 \geq N_2} \varepsilon^{1/2-\delta}\min(1,\varepsilon N_1^{-1}) \min(N_2^{1/2-\delta}, \eps N_2^{-1/2-\delta}) .
\end{equation*}
We break the right-hand side term into three sums
\begin{equation*}
\sum_{N_1 \geq N_2} = \sum_{N_1 \geq N_2 \geq \varepsilon} + \sum_{N_1 \geq \varepsilon \geq N_2} + \sum_{\varepsilon \geq N_1 \geq N_2}= I + II+III
\end{equation*}
which we bound as follows
\begin{align*}
I &= \sum_{N_1 \geq N_2 \geq \varepsilon} \varepsilon^{1/2-\delta} \varepsilon N_1^{-1}\varepsilon N_2^{-1/2-\delta}
   \lesssim \varepsilon^{1-2\delta}\\
II &= \sum_{N_1 \geq \varepsilon \geq N_2} \varepsilon^{1/2-\delta}\varepsilon N_1^{-1} N_2^{1/2-\delta}
   \lesssim \varepsilon^{1-2\delta}\\
III &= \sum_{\varepsilon \geq N_1 \geq N_2} \varepsilon^{1/2-\delta} N_2^{1/2-\delta}
   \lesssim \varepsilon^{1-2\delta}.
\end{align*}
The lemma follows.
\end{proof}

%%%%%%%%%%%%%%%%%%%%%%%%%%%%%%%%%%%%%%%%%%%%%%%%%%%%%%%%%%%%%%%%%%%%%%%%%%%%%%%%%%%%%%%%%%%
%
%
%                                   Section
%
%
%%%%%%%%%%%%%%%%%%%%%%%%%%%%%%%%%%%%%%%%%%%%%%%%%%%%%%%%%%%%%%%%%%%%%%%%%%%%%%%%%%%%%%%%%%%

\section{Perturbation Theory}
As mentioned in the introduction, the Cauchy problem for
\eqref{schrodinger equation} is locally well-posed in $\ho_x
(\R^4)$.  Indeed, this well-posedness extends to any interval
where one has uniform control of the $\ls$-norm (see
Lemma \ref{lemma persistence regularity} below).  This section
describes variants of the local well-posedness theory.  In
particular we will be interested in when we can perturb a solution
(or near-solution, see \eqref{short near solution} below) in the
energy norm when we can control the solution in $\ls$ and the
error in a dual Strichartz space.

As a first step, we consider the case where the near-solution, the
error, and the free evolution of the perturbation are small in
spacetime norms, but allowed to be large in the energy norm.

\begin{lemma}[Short-time perturbations]\label{lemma short time}
Let $I$ be a compact time interval, and let $\tilde u$ be a
near-solution to \eqref{schrodinger equation} on $I\times\R^4$ in the sense that
\begin{equation}\label{short near solution}
(i \partial_t + \Delta)\tilde u = |\tilde u|^2 \tilde u + e
\end{equation}
for some function $e$.  Suppose that we also have the energy bound
\begin{equation}\label{short energy bound}
\|\tilde u\|_{L^\infty_t \dot H^1_x(I\times\R^4)} \leq E
\end{equation}
for some $E>0$.  Let $t_0 \in I$ and let $u(t_0)$ close to $\tilde
u(t_0)$ in the sense that
\begin{equation}\label{short close energy}
\|u(t_0)-\tilde u(t_0)\|_{\dot H^1_x} \leq E'
\end{equation}
for some $E'>0$.
Assume also the smallness conditions
\begin{align}
\|\nabla \tilde u\|_{L^6_tL_x^{12/5}(\ir)} & \leq \delta \label{short small 1}\\
\|\nabla e^{i(t-t_0)\Delta}(u(t_0)-\tilde u(t_0))\|_{L^6_tL_x^{12/5}(\ir)} & \leq
\varepsilon \label{short small 2}\\
\|\nabla e\|_{L^2_t L^{4/3}_x(\ir)} & \leq \varepsilon \label{short
small 3}
\end{align}
for some $0 < \varepsilon, \delta < \varepsilon_0$, where $\eps_0 =\eps_0(E,E')$ is a small positive constant.

Then there exists a solution $u$ to \eqref{schrodinger equation}
on $\ir$ with the specified initial data $u(t_0)$ at $t_0$ satisfying
\begin{align}
\|u - \tilde u\|_{\lsi} & \lesssim \eps\label{short bound 3}\\
\|u - \tilde u \|_{\soi} & \lesssim \eps + E'\label{short bound 1}\\
\|u\|_{\soi} & \lesssim E' +E\label{short bound 2}\\
\bigl\|\nabla \bigl[(i \partial_t + \Delta)(u - \tilde u)+e\bigr]\bigr\|_{L^2_t L_x^{4/3} (\ir)}
& \lesssim \eps. \label{short bound 4}
\end{align}
\end{lemma}

\begin{remark}
Notice that $$\eqref{short small 2} \lesssim
\|e^{i(t-t_0)\Delta}(u(t_0)-\tilde u(t_0))\|_{\soi} \lesssim
\|u(t_0)-\tilde u(t_0)\|_{\dot H^1_x } \leq E',$$
so \eqref{short small 2} is redundant if $E' = O(\eps)$.
\end{remark}

\begin{proof}
By the well-posedness theory it will suffice to prove \eqref{short
bound 1} -- \eqref{short bound 4} as \textit{a priori} estimates,
that is, we will assume that the solution $u$ already exists and is
smooth on $I$.  Without loss of generality we may assume that $t
\geq t_0$, since the proof on the $t \leq t_0$ portion of $I$ is similar.

Let $v := u - \tilde u$, and for $t \in I$ define
\begin{align*}
S(t) := \bigl\|\nabla \bigl[(i \partial_t + \Delta)v + e\bigr]\bigr\|_{L^2_t L_x^{4/3} ([t_0,t]\times \R^4)}.
\end{align*}
By Lemma \ref{lemma linear strichartz}, \eqref{short small 2}, and \eqref{short small 3} we get
\begin{align}\label{Q bound}
\|\nabla v\|_{L_t^6L_x^{12/5} ([t_0, t]\times \R^4)}
&\lesssim \|\nabla e^{i(t-t_0)\Delta }v(t_0)\|_{L_t^6L_x^{12/5}([t_0,t]\times \R^4)}\\
&\quad+ \bigl\|\nabla\bigl[(i \partial_t +\Delta)v+ e\bigr]\bigr\|_{L^2_t L_x^{4/3}([t_0,t]\times \R^4)}\notag\\
&\quad+ \|\nabla e\|_{L^2_t L^{4/3}_x([t_0, t]\times \R^4)}\notag\\
&\lesssim S(t) + \eps.\notag
\end{align}
By Sobolev embedding, \eqref{Q bound} yields
\begin{align}\label{Q' bound}
\| v\|_{L_{t,x}^6([t_0, t]\times \R^4)}
\lesssim \|\nabla v\|_{L_t^6L_x^{12/5} ([t_0, t]\times \R^4)}
\lesssim S(t) + \eps.
\end{align}
On the other hand, since
$$(i \partial_t + \Delta)v = |\tilde u + v|^2 (\tilde u + v) -|\tilde u|^2 \tilde u - e = \sum_{j=1}^3 \O(v^j \tilde u^{3-j}) - e$$
we get
$$S(t) \lesssim \|\nabla \sum_{j=1}^3 \O (v^j \tilde u ^{3-j})\|_{L^2_t L_x^{4/3} ([t_0,t]\times \R^4)}.$$
Using the trilinear Strichartz estimate \eqref{trilinear 4/3} together with \eqref{short small 1},
\eqref{Q bound}, and \eqref{Q' bound}, one estimates
$$
S(t) \lesssim \eps + (E+ \eps+ \delta) S(t)^2 + S(t)^3.
$$
A standard continuity argument shows that if we take $\eps_0 = \eps_0(E,E')$ sufficiently small we obtain
\begin{align}\label{S bound}
S(t) \lesssim \eps
\end{align}
for all $t \in I$, which implies \eqref{short bound 4}.  Using \eqref{Q' bound} and \eqref{S bound},
one easily derives \eqref{short bound 3}. To obtain \eqref{short bound 1}, we use Strichartz's inequality:
\begin{align*}
\|u-\tilde{u}\|_{\dot{S}^1(\ir)}
&\lesssim \|u(t_0)-\tilde{u}(t_0)\|_{\dot{H}^1_x}+\bigl\|\nabla\bigl[(i \partial_t +\Delta)v+ e\bigr]\bigr\|_{L^2_t L_x^{4/3}(\ir)}\\
&\qquad+\|\nabla e\|_{L^2_t L_x^{4/3}(\ir)}\\
&\lesssim E'+\eps.
\end{align*}
By the triangle inequality, \eqref{short small 1} and \eqref{short bound 3} imply
$\|u\|_{L_{t,x}^6(\ir)}\lesssim \eps+\delta$. An application of Strichartz's inequality yields
\begin{align*}
\|u \|_{\soi} & \lesssim \| u (t_0)\|_{\dot H^1_x} + \|\nabla(| u|^2 u)\|_{L^2_t L^{4/3}_x(\ir)}\\
&\lesssim \| \tilde u (t_0)\|_{\dot H^1_x}+\| u (t_0)-\tilde u(t_0)\|_{\dot H^1_x}+\|u\|_{L_{t,x}^6(\ir)}^2\|u\|_{\soi}\\
& \lesssim E+E'+(\eps+\delta)^2\|u\|_{\soi},
\end{align*}
which proves \eqref{short bound 2}, provided $\eps_0$ is chosen sufficiently small
depending on $E$ and $E'$.
\end{proof}

We will also need the following version of the above lemma that
deals with near-solutions with large but finite $\ls$-norms.

\begin{lemma}[Long-time perturbations]\label{lemma long time}
Let $I$ be a compact time interval, and let $\tilde u$ be function
on $\ir$ that obeys the bounds
\begin{equation}\label{long l6 bound}
\|\tilde u\|_{L^6_{t,x}(\ir)} \leq M
\end{equation}
\begin{equation}\label{long energy bound}
\|\tilde u\|_{L^\infty_t \dot H^1_x(\ir)} \leq E
\end{equation}
for some $M,E>0$.  Suppose also that $\tilde u$ is a near-solution
to \eqref{schrodinger equation} in the sense of \eqref{short near
solution} for some function $e$. Let $t_0 \in I$ and let $u(t_0)$
close to $\tilde u(t_0)$ in the sense that
\begin{equation}\label{long close energy}
\|u(t_0)-\tilde u(t_0)\|_{\dot H^1_x} \leq E'
\end{equation}
for some $E'>0$.
Assume also the smallness conditions
\begin{align}
\|\nabla e^{i(t-t_0)\Delta}(u(t_0)-\tilde u(t_0))\|_{L^6_tL_x^{12/5}(\ir)} & \leq
\varepsilon \label{long small 1}\\
\|\nabla e\|_{L^2_t L^{4/3}_x(\ir)} & \leq \varepsilon \label{long
small 2}
\end{align}
for some $0 < \varepsilon < \varepsilon_1$, where $\eps_1 =\eps_1(E,E',M) $ is a small positive constant.

Then there exists a solution $u$ to \eqref{schrodinger equation}
on $\ir$ with the specified data $u(t_0)$ at $t_0$ satisfying
\begin{align}
\|u - \tilde u\|_{\lsi} & \leq C(E,E',M) \eps\label{long bound 3}\\
\|u - \tilde u \|_{\soi} & \leq C(E,E',M) \label{long bound 1}\\
\|u\|_{\soi} & \leq C(E,E',M)\label{long bound 2}.
\end{align}
\end{lemma}

\begin{remark}
The same computation as in Remark 3.1 shows that assumption \eqref{long small 1} is redundant if
one assumes $E' = O(\eps)$.  Also notice that if we take $e = 0$ the lemma implies local well-posedness in the energy
space whenever the $\ls$-norm is finite.
\end{remark}

\begin{proof}
Without loss of generality we may assume $t_0 = \inf I$.  Let $\eps_0 = \eps_0 (E, 2E')$ be as in the previous lemma
(we must replace $E'$ by $2E'$ because the kinetic energy of $u - \tilde u$ will grow in time).

We first note that \eqref{long l6 bound} implies $u\in \dot{S}^1(\ir)$. Indeed, subdividing $I$ into
$N_0 \sim (1 + \frac{M}{\eps_0})^{6}$ subintervals $J_k$ such
that on each $J_k$ we have
$$
\|\tilde u\|_{\ls (J_k \times \R^4)} \leq \eps_0,
$$
and using Lemma \ref{lemma linear strichartz} and Lemma \ref{lemma trilinear strichartz}, we estimate
\begin{align*}
\|\tilde u \|_{\so(J_k\times\R^4)} & \lesssim \|\tilde u (t_0)\|_{\dot H^1_x} +
\|\nabla(|\tilde u|^2 \tilde u)\|_{L^2_t L^{4/3}_x(J_k\times\R^4)} + \|\nabla
e\|_{L^2_t L^{4/3}_x(J_k\times\R^4)}\\
& \lesssim E + \|\tilde u\|_{\ls(J_k\times\R^4)}^2 \|\tilde u \|_{\so(J_k\times\R^4)} + \eps\\
&\lesssim E+\eps_0^2\|\tilde u \|_{\so(J_k\times\R^4)} + \eps.
\end{align*}
Thus, for $\eps_0$ sufficiently small, this yields
\begin{align*}
\|\tilde u \|_{\so(J_k\times\R^4)} \lesssim E + \eps.
\end{align*}
Summing these bounds over all the intervals $J_k$ we obtain
$$
\|\tilde u\|_{\soi}\leq C(M,E,\eps_0),
$$
which also implies by Lemma \ref{lemma strichartz norms} that
$$
\|\nabla \tilde u\|_{L_t^6L_x^{12/5}(\ir)}\leq C(M,E,\eps_0).
$$

We can now subdivide $I$ into $N_1=N_1(M, E, \eps_0)$ subintervals $I_j = [t_j , t_{j+1}]$ such
that on each $I_j$ we have
$$\|\tilde u\|_{\ls (I_j \times \R^4)} \leq \eps_0.$$
Choosing $\eps_1$ sufficiently small depending on $\eps_0,N_1,E, \text{ and } E',$ Lemma \ref{lemma short time} implies that
for each $j$ and all $0<\eps<\eps_1$,
\begin{align*}
\|u - \tilde u\|_{\ls (I_j \times \R^4)} & \leq C(j)\eps\\
\|u - \tilde u \|_{\so (I_j \times \R^4)} & \leq C(j)(\eps + E')\\
\|u\|_{\so (I_j \times \R^4)} & \leq C(j)( E' +E)\\
\bigl\|\nabla\bigl[(i \partial_t + \Delta)(u - \tilde u)+e\bigr]\bigr\|_{L^2_t L_x^{4/3} (I_j\times \R^4)} & \leq C(j)\eps,
\end{align*}
provided we can show that \eqref{long close energy} and \eqref{long small 1} hold when $t_0$ is replaced by $t_j$.
We check this using an inductive argument. By Strichartz we have the bounds
\begin{align*}
\|u(t_{j+1})-\tilde u(t_{j+1})\|_{\dot H^1_x}
&\lesssim \|u(t_0) -\tilde u(t_0)\|_{\dot H^1_x} + \| \nabla e \|_{L^2_t L^{4/3}_x(\ir)}\\
&\quad + \bigl\| \nabla \bigl[(i \partial_t + \Delta)(u - \tilde u) + e \bigr]\bigr\|_{L^2_t L_x^{4/3} ([t_0,t_{j+1}] \times\R^4)}\\
&\lesssim E'+ \sum_{k=0}^{j} C(k)\eps
\end{align*}
and similarly
\begin{align*}
\|\nabla e^{i(t-t_{j+1})\Delta}(&u(t_{j+1})-\tilde u(t_{j+1}))\|_{L_t^6L_x^{12/5}(\ir)} \\
&\lesssim \|\nabla e^{i(t-t_0)\Delta}(u(t_0)-\tilde u(t_0))\|_{L_t^6L_x^{12/5}(\ir)} + \| \nabla e \|_{L^2_t L^{4/3}_x(\ir)}\\
&\quad + \bigl\| \nabla \bigl[(i \partial_t + \Delta)(u - \tilde u) + e \bigr]\bigr\|_{L^2_t L_x^{4/3} ([t_0,t_{j+1}] \times\R^4)} \\
&\leq \sum_{k=0}^{j} C(k)\eps
\end{align*}
where $C(k)$ depends on $k, E, E', \text{ and }\eps_0$. Choosing $\eps_1$ sufficiently small depending on
$\eps_0,N_1,E, \text{ and } E',$ we can continue the inductive argument.
\end{proof}

\begin{remark}
The dependence of $\eps_1$ on the parameters $M,E, \text{ and } E'$ is extremely
bad.  Specifically, $\eps_1(M,E,E') \approx \exp (-M^C \langle E \rangle^C
\langle E' \rangle^C)$.  Since we'll use this lemma frequently, the bounds
in Theorem \ref{lemma THE THEOREM} will grow quite rapidly in $E$.  For remarks
on the final bound see Section 8.
It is of some interest to determine better bounds for the theorem,
but for us it will suffice that they remain finite.
\end{remark}

We'll also need the following related result concerning the
persistence of $L^2$, $\dot H_x^1$, and $\dot H_x^2$ regularity.

\begin{lemma}[Persistence of regularity]\label{lemma persistence regularity}
Let $k = 0,1,2$, $I$ a compact time interval, and $u$ a
finite-energy solution to \eqref{schrodinger equation} obeying
$$\|u\|_{\lsi} \leq M.$$
Then if $t_0 \in I$ and $u(t_0) \in \dot H_x^k$, we have
\begin{equation}\label{persistence regularity}
\|u\|_{\ski} \leq C(M,E(u))\|u(t_0)\|_{\dot H_x^k}.
\end{equation}
\end{lemma}

So once we have $\ls$ control of a finite-energy solution, we
control all Strichartz norms as well.  If the initial data is
$\htwo_x$ then we also control the $\st$-norm.  By iterating, we can
continue a local-in-time solution to any interval on which we have
uniform control of the $\ls$-norm.

\begin{proof}
Again it suffices to prove \eqref{persistence regularity} as an
\textit{a priori} bound.  Let $\tilde u = u$, $e = 0$, and $E' =
0$ and apply Lemma \ref{lemma long time} to get
$$\|u\|_{\soi} \lesssim C(M,E).$$
By \eqref{trilinear strichartz} we also have
\begin{equation}\label{grad k 12 bound}
\|\nabla^k \O(u^3)\|_{L^1_t L^2_x(I\times\R^4)} \lesssim \|u\|_{\ls(I\times\R^4)}
\|u\|_{\sk(I\times\R^4)} \|u\|_{\so(I\times\R^4)}.
\end{equation}
Now, divide $I$ into $N \approx (1 + \frac{M}{\delta})^{6}$
subintervals $I_j = [T_j , T_{j+1}]$ on which
$$\|u\|_{\ls (I_j \times \R^4)} \leq \delta$$
where $\delta$ will be chosen later.  On each $I_j$
Strichartz and \eqref{grad k 12 bound} yield
\begin{align*}
\|u\|_{\sk (I_j \times \R^4)} &\lesssim \|u(T_j)\|_{\dot
H_x^k(\R^4)} + \|\nabla^k(|u|^2 u)\|_{L^1_t L^2_x (I_j \times
\R^4)}\\
&\lesssim \|u(T_j)\|_{\dot H_x^k(\R^4)} + \|u\|_{\ls (I_j \times
\R^4)} \|u\|_{\sk (\ir)} \|u\|_{\so (\ir)}.
\end{align*}
So choosing $\delta \leq \frac{1}{2} C(M,E)^{-1}$ we get
\begin{equation}\label{sk hk bound}
\|u\|_{\sk (I_j \times \R^4)} \lesssim \|u(T_j)\|_{\dot H_x^k
(\R^4)}.
\end{equation}
The lemma follows from adding up the bounds \eqref{sk hk bound} on
each subinterval.
\end{proof}

%%%%%%%%%%%%%%%%%%%%%%%%%%%%%%%%%%%%%%%%%%%%%%%%%%%%%%%%%%%%%%%%%%%%%%%%%%%%%%%%%%%%%%%%%%%
%
%
%                                   Section
%
%
%%%%%%%%%%%%%%%%%%%%%%%%%%%%%%%%%%%%%%%%%%%%%%%%%%%%%%%%%%%%%%%%%%%%%%%%%%%%%%%%%%%%%%%%%%%

\section{Frequency Localization and Space Concentration}

Recall from the introduction that we expect a minimal energy
blowup solution to be localized in both frequency and space.  In
this section we will prove that this is indeed the case (we will
not actually prove that the solution is localized in space, just
that it concentrates; see the discussion after the proof of
Corollary \ref{lemma freq loc}). The first step is the following
proposition:

\begin{proposition}[Frequency delocalization $\Rightarrow$ spacetime
bound]\label{lemma freq loc imp spacetime bds} Let $u$ be a solution to \eqref{schrodinger equation} on $I_*\times\R^4$ with $E(u)\leq E_{crit}$.
Let $\eta > 0$ and suppose there exist a dyadic frequency $N_{lo}>0$ and a time $t_0
\in I_*$ such that we have the energy separation conditions
\begin{equation}\label{freq deloc lo freq}
\|P_{\leq N_{lo}} u(t_0)\|_{\ho_x} \geq \eta
\end{equation}
and
\begin{equation}\label{freq deloc hi freq}
\|P_{\geq K(\eta)N_{lo}} u(t_0)\|_{\ho_x} \geq \eta.
\end{equation}
If $K(\eta)$ is sufficiently large depending on $\eta$ we have
\begin{equation}\label{freq deloc ls}
\|u\|_{\ls (I_* \times \R^4)}\leq C(\eta).
\end{equation}
\end{proposition}

\begin{proof}
Let $0 < \eps = \eps(\eta) \ll 1$ be a small number to be chosen later.  If $K(\eta)$ is sufficiently large depending
on $\eps$ ($K(\eta)$ needs to be of the order $\eps^{-(\eps^{-2})}$), then one can find $\eps^{-2}$ disjoint intervals $[\eps^2 N_j,
\eps^{-2} N_j]$ contained in $[N_{lo}, K(\eta)N_{lo}]$.  By \eqref{kinetic energy bound}
and the pigeonhole principle, we may find an $N_j$ such that the interval $[\eps^2 N_j, \eps^{-2}N_j]$ has very little energy
$$\|P_{\eps^2 N_j \leq \cdot \leq \eps^{-2}N_j}u(t_0)\|_{\ho_x} \lesssim \eps.$$ Since both the statement and conclusion of
the proposition are invariant under the scaling \eqref{scaling}, we normalize $N_j = 1$.

Define $\ulo(t_0) := P_{\leq \eps}u(t_0)$ and $\uhi(t_0) = P_{\geq
\eps^{-1}}u(t_0)$.  We claim that $\uhi$ and $\ulo$ have smaller
energy than $u$.

\begin{lemma}\label{lemma uhi ulo smaller energy}
If $\eps$ is sufficiently small depending on $\eta$, then we have
$$E(\ulo(t_0)),E(\uhi(t_0)) \leq E_{crit} - c\eta^C.$$
\end{lemma}
\begin{proof}
We'll prove this for $\ulo$; the proof for $\uhi$ is similar.
Define $\uhip(t_0) := P_{> \eps}u(t_0)$ so that $u(t_0) = \ulo
(t_0) + \uhip(t_0)$, and consider the quantity
\begin{equation}\label{freq deloc energy 1}
|E(u(t_0)) - E(\ulo (t_0)) - E(\uhip(t_0))|.
\end{equation}
By the definition of energy we can bound this by
\begin{equation}\label{freq deloc energy 2}
|\langle \nabla \ulo (t_0),\nabla \uhip (t_0) \rangle| +
|\int(|u(t_0)|^4-|\ulo(t_0)|^4-|\uhip(t_0)|^4) dx|.
\end{equation}
We deal with the potential energy term first.  By the pointwise estimate
$$||u(t_0)|^4-|\ulo(t_0)|^4-|\uhip(t_0)|^4| \lesssim
|\ulo (t_0)||\uhip (t_0)|(|\ulo (t_0)|+|\uhip (t_0)|)^2$$ and
H\"older, we can bound the potential energy term in \eqref{freq deloc
energy 2} by
\begin{equation}\label{freq deloc energy 3}
\|\ulo (t_0)\|_{L^\infty_x}\|\uhip (t_0)\|_{L^2_x} (\|\ulo
(t_0)\|_{L^4_x} + \|\uhip (t_0)\|_{L^4_x})^2.
\end{equation}
An application of Bernstein yields
\begin{align*}
\|\ulo (t_0)\|_{L^\infty_x}\lesssim \eps \|\ulo(t_0)\|_{L_x^4}\lesssim \eps\|\ulo\|_{\dot{H}^1_x}\lesssim \eps.
\end{align*}
Similarly
\begin{align*}
\|\uhip (t_0)\|_{L^2_x}
   &\lesssim \sum_{N > \eps} \|P_N u(t_0)\|_{L^2_x}
     \lesssim \sum_{N > \eps}N^{-1} \|P_N u (t_0)\|_{\ho_x}\\
   &\lesssim \sum_{N > 1/\eps^2} N^{-1} + \sum_{\eps < N \leq 1/\eps^2}N^{-1}\eps
     \lesssim \eps^2 + \eps(1/\eps - \eps^2)\\
   &\lesssim 1.
\end{align*}
Thus $\eqref{freq deloc energy 3} \lesssim \eps$.

Now we deal with the kinetic part of \eqref{freq deloc energy 2}. We estimate
\begin{align*}
|\langle \nabla \ulo (t_0),\nabla \uhip (t_0) \rangle|
   &\lesssim |\langle \nabla P_{>\eps}P_{\leq\eps}u(t_{0}), \nabla u(t_0)\rangle| \\
   &\lesssim \|\nabla P_{>\eps}P_{\leq\eps}u(t_{0})\|_{L^2_x} \|\nabla u(t_0)\|_{L^2_x}.
\end{align*}
As
$$
\|\nabla u(t_0)\|_{L^2_x}\lesssim 1
$$
and
\begin{align*}
\|\nabla P_{>\eps} P_{\leq\eps}u(t_{0})\|_{L^2_x}
   &=\|(\nabla P_{>\eps}P_{\leq\eps}u(t_{0}))^{\wedge}\|_{L^2_x} \\
   &=\|\xi \varphi(\xi/\eps)(1-\varphi(\xi/\eps))\widehat{u(t_0)}(\xi)\|_{L^2_x} \\
   &\lesssim \eps \|\widehat{u_{hi'}(t_0)}\|_{L^2_x}
     \lesssim \eps
\end{align*}
we get as a final bound on kinetic energy
$$
|\langle \nabla \ulo (t_0),\nabla \uhip (t_0) \rangle|\lesssim \eps.
$$

Therefore $\eqref{freq deloc energy 1} \lesssim \eps$. As $$E(u) \leq E_{crit},$$ and by
hypothesis $$E(\uhip(t_0))\gtrsim \|\uhip(t_0)\|^2_{\ho_x} \gtrsim \eta^2,$$ the triangle inequality implies
$E(\ulo(t_0)) \leq E_{crit} -c\eta^C$.
\end{proof}

Now since $E(\ulo(t_0)),E(\uhi(t_0)) \leq \ecrit - c\eta^C < \ecrit$
we can apply Lemma
\ref{lemma induct on energy} to deduce that there exist solutions
$\ulo$ and $\uhi$ on the slab $\is$ with initial data $\ulo (t_0)$
and $\uhi (t_0)$ such that
\begin{align*}
\|\ulo\|_{\sois} &\lesssim C(\eta)\\
\|\uhi\|_{\sois} &\lesssim C(\eta).
\end{align*}
Define $\tilde u := \ulo + \uhi$.  We claim that $\tilde u$ is a
near-solution to \eqref{schrodinger equation}.

\begin{lemma}\label{lemma uhi plus ulo near soln}
We have
$$i \tilde u_t + \Delta \tilde u = |\tilde u|^2 \tilde u - e$$
where the error $e$ obeys the bound
\begin{equation}\label{freq deloc error bound}
\|\nabla e\|_{L^2_t L^{4/3}_x (\is)} \lesssim C(\eta)\eps^{1/2}.
\end{equation}
\end{lemma}

\begin{proof}
By the above estimates and \eqref{mass high freq bound} we have
that $\|\uhi(t_0)\|_{L^2_x} \lesssim \eps$, and so by \eqref{persistence regularity} we get
\begin{equation}
\|\uhi\|_{\szis} \lesssim C(\eta)\eps.
\end{equation}
Similarly, from \eqref{kinetic energy bound} and  Bernstein  we have
$$\|\ulo(t_0)\|_{\dot{H}^2_x} \lesssim \eps \|\ulo (t_0)\|_{\dot{H}^1_x}
\lesssim \eps \|\ulo \|_{L^{\infty}_t\dot{H}^1_x(I_*\times\R^4)}
\lesssim C(\eta)\eps,$$ and so by \eqref{persistence regularity}
\begin{equation}
\|\ulo\|_{\stis} \lesssim C(\eta)\eps.
\end{equation}
From Lemma \ref{lemma trilinear strichartz} we have the additional
bounds
\begin{align*}
\| |\uhi|^2\uhi\|_{L^1_t L^2_x (\is)} &\lesssim C(\eta)\eps\\
\| \nabla(|\uhi|^2\uhi)\|_{L^1_t L^2_x (\is)} &\lesssim C(\eta)\\
\| \nabla(|\ulo|^2\ulo)\|_{L^1_t L^2_x (\is)} &\lesssim C(\eta)\\
\| \nabla^2(|\ulo|^2\ulo)\|_{L^1_t L^2_x (\is)} &\lesssim
C(\eta)\eps.
\end{align*}
Applying Lemma \ref{lemma trilinear strichartz hilo} we see
$$\|\nabla \O(\uhi^j \ulo^{3-j})\|_{L^2_t L^{4/3}_x (\is)}
\lesssim C(\eta)\eps^{1/2}$$ for $j = 1,2$.  Since $e =
\sum_{j=1}^2 \O(\uhi^j \ulo^{3-j})$ the claim follows.
\end{proof}

We derive estimates on $u$ from those on $\tilde u$ via
perturbation theory.  More precisely, we know that
$$\|u(t_0)-\tilde u(t_0)\|_{\ho_x} \lesssim \eps$$ and
$$\|\tilde u \|_{\lsis} \lesssim \|\ulo\|_{\so(I_* \times\R^4) } + \|\uhi\|_{\so(I_* \times\R^4)}
\lesssim C(\eta).$$ So if $\eps$ is sufficiently small depending
on $\eta$, we can apply Lemma \ref{lemma long time} to deduce the
bound \eqref{freq deloc ls}.  This concludes the proof of
Proposition \ref{lemma freq loc imp spacetime bds}.
\end{proof}

Comparing \eqref{freq deloc ls} with \eqref{l6 HUGE} gives the
desired contradiction if $u$ satisfies the hypotheses of
Proposition \ref{lemma freq loc imp spacetime bds}.  We therefore
expect $u$ to be localized in frequency for each $t$. Indeed we
have:

\begin{corollary}[Frequency localization of energy at each time]\label{lemma freq loc}
Let $u$ be a minimal energy blowup solution of \eqref{schrodinger
equation}.  Then for each time $t_0 \in I_*$ there exists a dyadic
frequency $N(t)\in 2^{\Z}$ such that for every $\eta_3 \leq \eta
\leq \eta_0$ we have small energy at frequencies $\ll N(t)$
\begin{equation}
\|P_{\leq c(\eta)N(t)}u(t)\|_{\ho_x} \leq \eta
\end{equation}
small energy at frequencies $\gg N(t)$
\begin{equation}
\|P_{\geq C(\eta)N(t)}u(t)\|_{\ho_x} \leq \eta
\end{equation}
and large energy at frequencies $\sim N(t)$
\begin{equation}
\|P_{c(\eta)N(t)<\cdot<C(\eta)N(t)}u(t)\|_{\ho_x} \sim 1
\end{equation}
where the values of $0 < c(\eta) \ll 1 \ll C(\eta) < \infty$
depend on $\eta$.
\end{corollary}

\begin{proof}
For $t \in I_*$ define $$N(t) := \sup \{N \in 2^{\Z} : \|P_{\leq
N} u(t) \|_{\ho_x} \leq \eta_0 \}.$$  As $u$ is Schwartz, $N(t) >
0$; as $\|u\|_{L^\infty_t \ho_x} \sim 1$, $N(t) < \infty$. From
the definition of $N(t)$ we have that $$\|P_{\leq 2N(t)}
u(t)\|_{\ho_x} > \eta_0.$$  Let $\eta_3 \leq \eta \leq \eta_0$. If
$C(\eta) \gg 1$ then we must have $\| P_{\geq C(\eta)N(t)}
u(t)\|_{\ho_x} \leq \eta$, since otherwise Proposition \ref{lemma
freq loc imp spacetime bds} would imply $\|u\|_{\lsis} \lesssim
C(\eta)$, which would contradict $u$ being a minimal energy blowup
solution.

Similarly, $\|u\|_{L^\infty_t \ho_x} \sim 1$ implies that
$$\|P_{c(\eta_0)N(t) < \cdot < C(\eta_0)N(t)} u(t_0) \|_{\ho_x}
\sim 1$$ and therefore that $$\|P_{c(\eta)N(t)< \cdot <
C(\eta)N(t)} u(t_0)\|_{\ho_x} \sim 1$$ for all $\eta_3 \leq \eta
\leq \eta_0$.  Thus, if $c(\eta) \ll 1$ then $\|P_{\leq
c(\eta)N(t)} u(t)\|_{\ho_x} \leq \eta$ for all $\eta_3 \leq \eta
\leq \eta_0$, since otherwise Proposition \ref{lemma freq loc imp
spacetime bds} would again imply $\|u\|_{\lsis} \lesssim C(\eta)$.
\end{proof}

Having shown that a minimal energy blowup solution must be
localized in frequency, we turn our attention to space.  While it
is true that such a solution is localized in space (the methods employed in \cite{ckstt:gwp}
carry through to the four dimensional
energy-critical NLS case), we will only need the following weaker
result concerning concentration of the solution (roughly,
\emph{concentration} will mean large at some point, while we
reserve \emph{localization} to mean simultaneously concentrated
and small at points far from the concentration point).  To obtain the concentration result, we divide the
interval $I_*$ into three consecutive subintervals $I_* = I_{-}
\cup I_0 \cup I_{+}$, each containing a third of the $\ls$ density
of $u$:
$$\int_I \int_{\R^4} |u(t,x)|^6 dx dt = \frac{1}{3} \int_{I_*}
\int_{\R^4} |u(t,x)|^6 dx dt \text{ for } I = I_-,I_0,I_+.$$ It is
on the middle interval $I_0$ that we will show physical space
concentration.  The first step is:

\begin{proposition}[Potential energy bounded from below]\label{lemma potential bdd below}
For any minimal energy blowup solution of \eqref{schrodinger
equation} and all $t \in I_0$ we have
\begin{equation}\label{potential bdd below}
\|u(t)\|_{L^4_x} \geq \eta_1.
\end{equation}
\end{proposition}

\begin{proof}
We will use an idea of Bourgain \cite{borg:scatter}.  If the linear
evolution of the solution does not concentrate at some point in
spacetime, then we can use the small data theory and iterate. So
say the linear evolution concentrates at some point $(t_1, x_1)$.  If the
linear evolution is small in $L^4_x$ at time $t_0$, we show $t_0$
must be far from $t_1$.  We then remove the energy concentrating
at $(t_1, x_1)$ and use induction on energy.

More formally, we'll argue by contradiction.  Suppose there exists some time $t_0 \in I_0$ such that
\begin{equation}\label{contra for pot bdd below}
\|u(t_0)\|_{L^4_x} \lesssim \eta_1.
\end{equation}
Using \eqref{scaling} we scale $N(t_0)=1$.  If the
linear evolution $e^{i(t-t_0)\Delta}u(t_0)$ had small $L^6_{t,x}$-norm,
then by perturbation theory the nonlinear solution
would have small $L^6_{t,x}$-norm as well.  Hence, we may assume
$\|e^{i(t-t_0)\Delta}u(t_0)\|_{L^6_{t,x}(\rr)} \gtrsim 1$.

On the other hand Corollary \ref{lemma freq loc} implies that
$$\| \pl u(t_0) \|_{\ho_x} + \| \ph u(t_0)\|_{\ho_x} \lesssim \eta_0,$$
where $\pl = P_{<c(\eta_0)}$ and $\ph = P_{>C(\eta_0)}$.  Then by
Strichartz
$$\|\propagate \pl u(t_0)\|_{\lsr} + \|\propagate \ph
u(t_0)\|_{\lsr}\lesssim \eta_0.$$ Thus $$\|\propagate \pmed
u(t_0)\|_{\lsr} \sim 1,$$ where $\pmed = 1 - \pl - \ph$.  However, $\pmed u(t_0)$ has bounded energy (by \eqref{kinetic
energy bound}) and Fourier support in $c(\eta_0) \lesssim |\xi|
\lesssim C(\eta_0)$.  Another application of Strichartz yields
$$\|\propagate \pmed u(t_0)\|_{L^3_{t,x}} \lesssim \|\pmed
u(t_0)\|_{L^2_x} \lesssim C(\eta_0).$$  Combining these estimates
with H\"older we have $$\|\propagate \pmed
u(t_0)\|_{L^\infty_{t,x}} \gtrsim c(\eta_0).$$  In particular
there exist a time $t_1 \in \R$ and a point $x_1 \in \R^4$ so that
$$|e^{i(t_1-t_0)\Delta} (\pmed u(t_0))(x_1)| \gtrsim c(\eta_0).$$
We may perturb $t_1$ so that $t_1 \neq t_0$, and by time reversal
symmetry we may take $t_1 < t_0$.  Let $\delta_{x_1}$ be the Dirac
mass at $x_1$.  Define $f(t_1) := \pmed \delta_{x_1}$ and for $t > t_1$ define
$f(t) := \propagateo f(t_1)$.  One should think of $f(t_1)$ as basically $u$ at $(t_1,x_1)$.
The point is then to compare $u(t_0)$
to the linear evolution of $f(t_1)$ at time $t_0$. We will show that $f(t)$ decays rapidly
in any $L^p_x$-norm for $1\leq p\leq \infty$.

\begin{lemma}\label{lemma potential bdd below lemma}
For any $t \in \R$ and any $1 \leq p \leq \infty$ we have
$$\|f(t)\|_{L^p_x} \lesssim C(\eta_0) \langle t-t_1\rangle^{\frac{4}{p}-2}.$$
\end{lemma}

\begin{proof}
We translate so that $t_1 = x_1 = 0$, then use Bernstein and
the unitarity of $e^{it\Delta}$ to get
$$\|f(t)\|_{L^\infty_x} \lesssim C(\eta_0) \|f(t)\|_{L^2_x} =
C(\eta_0) \|\pmed \delta_{x_1}\|_{L^2_x} \lesssim C(\eta_0).$$
By \eqref{dispersive ineq} we also have $$\|f(t)\|_{L^\infty_x}
\lesssim |t|^{-2}\|\pmed \delta_{x_1}\|_{L^1_x} \lesssim C(\eta_0)
|t|^{-2}.$$  Combining these two we obtain
$$\|f(t)\|_{L^\infty_x} \lesssim C(\eta_0) \langle t\rangle^{-2}.$$  This proves the lemma in the case $p = \infty$.

For other $p$'s we use \eqref{fourier rep} to write $$f(t,x) =
\int e^{2 \pi i (x \xi - 2 \pi t |\xi|^2)} \phi_{med} (\xi) d\xi$$
where $\phi_{med}$ is the Fourier multiplier corresponding to
$\pmed$.  For $|x| \gg 1 + |t|$, repeated integration by parts
shows $|f(t,x)| \lesssim |x|^{-100}$.  On $|x| \lesssim 1 + |t|$,
one integrates using the above $L_x^\infty$-bound.
\end{proof}

{}From \eqref{contra for pot bdd below} and H\"older we have
$$| \langle u(t_0),f(t_0)\rangle| \lesssim
\|f(t_0)\|_{L^{4/3}_x}\|u(t_0)\|_{L^4_x} \lesssim \eta_1 C(\eta_0)
\langle t_1 - t_0\rangle.$$ On the other hand
$$| \langle u(t_0),f(t_0) \rangle | = | \langle e^{i(t_1 -
t_0)\Delta}\pmed u(t_0),\delta_{x_1}\rangle| \gtrsim c(\eta_0).$$
So since $\langle t_1 - t_0\rangle \gtrsim c(\eta_0)/\eta_1$, we have that
$t_1$ is far from $t_0$.  In particular, the time of concentration must be far from where the $L^4_x$-norm
is small.  Also, from Lemma \ref{lemma potential bdd below
lemma} we have that $\nabla f$ has a small $L^6_tL_x^{12/5}$-norm to the future
of $t_0$ (recall $t_1 < t_0$):
\begin{equation}\label{small future l6}
\|f\|_{L^6_tL_x^{12/5}([t_0,\infty)\times \R^4)} \lesssim C(\eta_0)|t_1- t_0|^{-1/6} \lesssim C(\eta_0)\eta_1^{1/6}.
\end{equation}
Now we use the induction hypothesis.  Split $u(t_0) = v(t_0) + w(t_0)$ where
$w(t_0) = \delta e^{i \theta} \Delta^{-1} f(t_0)$ for some small $\delta = \delta(\eta_0) > 0$ and phase $\theta$ to
be chosen later.  The point of $\Delta^{-1}$ in the definition of $w(t_0)$ is that our inner products will be in
$\dot{H}^1_x$ instead of the more common $L^2_x$. One should think of $w(t_0)$ as the contribution coming from the point $(t_1,x_1)$
where the solution concentrates. We will show that for an appropriate choice of
$\delta$ and $\theta$, $v(t_0)$ has slightly smaller energy than
$u$.  By the definition of $f$ and an integration by parts we have
\begin{align*}
\frac{1}{2} \int_{\R^4} |\nabla v(t_0)|^2 dx &= \frac{1}{2} \int
|\nabla u(t_0) - \nabla w(t_0)|^2 dx\\
&=\frac{1}{2} \int |\nabla u(t_0)|^2 dx - \delta \text{Re} \int
e^{-i \theta} \overline{\nabla \Delta^{-1}f(t_0)} \cdot \nabla
u(t_0)
dx\\
&\quad +O(\delta^2\|\Delta^{-1}f(t_0)\|_{\ho_x}^2)\\
& \leq \ecrit + \delta \text{Re } e^{-i\theta}\langle
u(t_0),f(t_0) \rangle + O(\delta^2 C(\eta_0)).
\end{align*}
Choosing $\delta$ and $\theta$ appropriately we get
$$\frac{1}{2} \int_{\R^4} |\nabla v(t_0)|^2 dx \leq \ecrit -
c(\eta_0).$$

Again by Lemma \ref{lemma potential bdd below lemma} we have
$$\|w(t_0)\|_{L^4_x} \lesssim C(\eta_0)\|f(t_0)\|_{L^4_x} \lesssim
C(\eta_0)\langle t_1 - t_0\rangle^{-1} \lesssim C(\eta_0)\eta_1.$$  So by
\eqref{contra for pot bdd below} and the triangle inequality we
obtain
$$\int_{\R^4} |v(t_0)|^4 dx \lesssim C(\eta_0)\eta_1^4.$$

Combining the above two energy estimates we see that $$E(v(t_0))
\leq \ecrit - c(\eta_0),$$ and Lemma \ref{lemma induct on energy} implies
there exists a global solution $v$ of \eqref{schrodinger equation} on
$[t_0, \infty) \times \R^4$ with data $v(t_0)$ at time $t_0$
satisfying $$\|v\|_{L^6_{t,x}([t_0, \infty) \times \R^4)} \leq
M(\ecrit - c(\eta_0)) = C(\eta_0).$$

However, \eqref{small future l6} and frequency localization give
$$
\|\nabla \propagate w(t_0)\|_{L^6_tL_x^{12/5}([t_0,\infty)\times \R^4)}
\lesssim \|\nabla \Delta^{-1} f\|_{L^6_tL_x^{12/5}([t_0,\infty)\times \R^4)}
\lesssim C(\eta_0) \eta_1^{1/6}.
$$
So if $\eta_1$ is sufficiently small depending on $\eta_0$, we can apply Lemma \ref{lemma long time} with
$\tilde u = v$ and $e = 0$ to conclude that $u$ extends to all of $[t_0, \infty)$ and obeys
$$
\|u\|_{L^6_{t,x}([t_0,\infty)\times\R^4)} \lesssim C(\eta_0,\eta_1).
$$
As $[t_0,\infty)$ contains $I_+$, the above estimate contradicts \eqref{l6 HUGE} if $\eta_4$ is chosen sufficiently small.
This concludes the proof of Proposition \ref{lemma potential bdd below}.
\end{proof}

Using \eqref{potential bdd below} we can deduce the desired concentration result:

\begin{proposition}[Spatial concentration of energy at each time]\label{lemma physical concentration}
For any minimal energy blowup solution of \eqref{schrodinger
equation} and for each $t \in I_0$, there exists $x(t) \in \R^4$
such that
\begin{equation}\label{physical conc kinetic}
\int_{|x-x(t)| \leq C(\eta_1)/N(t)} |\nabla u(t,x)|^2 dx \gtrsim
c(\eta_1)
\end{equation}
and
\begin{equation}\label{physical conc lp}
\int_{|x-x(t)| \leq C(\eta_1)/N(t)} |u(t,x)|^p dx \gtrsim
c(\eta_1) N(t)^{p-4}
\end{equation}
for all $1 < p < \infty$, where the implicit constants depend on
$p$. In particular
\begin{equation}\label{physical conc potential}
\int_{|x-x(t)| \leq C(\eta_1)/N(t)} |u(t,x)|^4 dx \gtrsim
c(\eta_1).
\end{equation}
\end{proposition}

Note that $u(t,x)$ is roughly of size $N(t)$ on average when $|x-x(t)|\leq N(t)^{-1}$,
which is consistent with both the concentration of energy and the uncertainty principle.

\begin{proof}
Fix $t$ and normalize $N(t)=1$.  By Corollary \ref{lemma freq
loc} we have
$$\|P_{<c(\eta_1)} u(t)\|_{\ho_x} + \|P_{>C(\eta_1)} u(t)\|_{\ho_x} \lesssim \eta_1^{100}  .$$
Sobolev embedding implies
$$\|P_{<c(\eta_1)} u(t)\|_{L^4_x} +\|P_{>C(\eta_1)} u(t)\|_{L^4_x}
\lesssim \eta_1^{100}$$ and so by \eqref{potential bdd below}
$$\|\pmed u(t)\|_{L^4_x} \gtrsim \eta_1,$$ where $\pmed = P_{c(\eta_1) \leq \cdot
\leq C(\eta_1)}$.  On the other hand, by \eqref{kinetic energy
bound} we have $$\|\pmed u(t)\|_{L^2_x} \lesssim C(\eta_1).$$
By H\"older $$\eta_1 \lesssim \| \pmed u(t)\|_{L^4_x} \lesssim
\|\pmed u(t)\|_{L^\infty_x}^{1/2} \cdot \| \pmed
u(t)\|_{L^2_x}^{1/2} \lesssim C(\eta_1) \|\pmed
u(t)\|_{L^\infty_x}^{1/2},$$ which implies that $$\|\pmed
u(t)\|_{L^\infty_x} \gtrsim c(\eta_1).$$ In particular, there
exists a point $x(t) \in \R^4$ so that
\begin{equation}\label{pmed bdd below}
c(\eta_1) \lesssim |\pmed u(t,x(t))|.
\end{equation}
As our function is now localized both in frequency and in space, all the Sobolev norms are practically equivalent.
So let's consider the operator $\pmed \nabla \Delta^{-1}$ and let $K_{med}$ denote its kernel.  Then
\begin{align*}
c(\eta_1) &\lesssim |\pmed u(t,x(t))| \lesssim |K_{med} * \nabla
u(t,x(t))| \lesssim \int |K_{med}(x(t)-x)||\nabla u(t,x)| dx\\
& \sim \int_{|x-x(t)|<C(\eta_1)}|K_{med}(x(t)-x)||\nabla u(t,x)|
dx\\
&\quad + \int_{|x-x(t)|\geq C(\eta_1)}|K_{med}(x(t)-x)||\nabla
u(t,x)|
dx\\
& \lesssim C(\eta_1) \Bigl( \int_{|x-x(t)|< C(\eta_1)}|\nabla
u(t,x)|^2 dx \Bigr)^{1/2} + \int_{|x-x(t)| \geq C(\eta_1)}
\frac{|\nabla u(t,x)|}{|x-x(t)|^{100}} dx ,
\end{align*}
where in order to obtain the last inequality we used Cauchy-Schwarz and that $K_{med}$ is a Schwartz
function.  Therefore, by \eqref{kinetic
energy bound} and possibly making $C(\eta_1)$ larger, we have
$$c(\eta_1) \lesssim \Bigl( \int_{|x-x(t)|< C(\eta_1)}|\nabla
u(t,x)|^2 dx \Bigr)^{1/2} + C(\eta_1)^{-\alpha}$$ for some $\alpha
> 0$, proving \eqref{physical conc kinetic}.

Now let $\tilde K_{med}$ be the kernel associated to $\pmed$, and
let $1 < p < \infty$.  As above we get
\begin{align*}
c(\eta_1)
&\lesssim \int |\tilde K_{med}(x(t)-x)|| u(t,x)| dx\\
& \sim \int_{|x-x(t)|<C(\eta_1)}|\tilde K_{med}(x(t)-x)||u(t,x)|dx\\
&\quad + \int_{|x-x(t)|\geq C(\eta_1)}|\tilde K_{med}(x(t)-x)|| u(t,x)|dx\\
& \lesssim C(\eta_1) \Bigl( \int_{|x-x(t)|< C(\eta_1)}|u(t,x)|^p dx \Bigr)^{1/p} \\
&\quad + \|u(t)\|_{L^4_x} \Bigl( \int_{|x-x(t)| \geq C(\eta_1)}\frac{1}{|x-x(t)|^{100 \cdot 4/3}} dx \Bigr)^{3/4}\\
& \lesssim C(\eta_1) \Bigl( \int_{|x-x(t)|< C(\eta_1)}|
u(t,x)|^p dx \Bigr)^{1/p} + C(\eta_1)^{-\alpha}
\end{align*}
for some $\alpha > 0$ which, after undoing the scaling, proves
\eqref{physical conc lp} if $C(\eta_1)$ is sufficiently large.
\end{proof}

%%%%%%%%%%%%%%%%%%%%%%%%%%%%%%%%%%%%%%%%%%%%%%%%%%%%%%%%%%%%%%%%%%%%%%%%%%%%%%%%%%%%%%%%%%%
%
%
%                                   Section
%
%
%%%%%%%%%%%%%%%%%%%%%%%%%%%%%%%%%%%%%%%%%%%%%%%%%%%%%%%%%%%%%%%%%%%%%%%%%%%%%%%%%%%%%%%%%%%

\section{Frequency-Localized Interaction Morawetz Inequality}

The goal of this section is to prove

\begin{proposition}[Frequency-localized interaction Morawetz
estimate]Roughly speaking, this proposition states that after throwing away some low frequency components of the minimal energy
blowup solution, the remainder obeys good $L^2_t\dot{H}^{-1/2}_x$ estimates. \label{prop flim}
Assuming u is a minimal energy blowup solution of
\eqref{schrodinger equation} and $N_{*}<c(\eta_{1})N_{min}$, we have
\begin{align}
\int_{I_{0}}\int_{\R^{4}}\int_{\R^{4}}\frac{|P_{\geq N_{*}}u(t,x)|^{2}|P_{\geq N_{*}}u(t,y)|^{2}}{|x-y|^{3}}dxdydt
\lesssim\eta_{1}N_{*}^{-3}.\label{flim}
\end{align}
\end{proposition}

In order to prove the proposition we introduce an interaction potential generalization of the classical Morawetz inequality.

\subsection{An interaction virial identity and a general interaction
Morawetz estimate for general equations}

We start by recalling the standard Morawetz action centered at a point for general equations. Let $a$  be
a function on the slab $I\times \R^{n}$ and $\phi$ satisfy $i\phi_{t}+\Delta \phi=\mathcal{N}$ on $I\times \R^{n}$.
We define the virial potential to be
$$V_{a}(t)=\int_{\R^{n}}a(x)|\phi(t,x)|^{2}dx$$
and the Morawetz action centered at zero to be
$$M_{a}^0(t)=2\int_{\R^{n}}a_{j}(x)\Im(\overline{\phi(x)}\phi_{j}(x))dx.$$
A computation shows that
$$\partial_{t}V_{a}=M_{a}^0+2\int_{\R^n} a\{\mathcal{N},\phi\}_{m}dx,$$
where the mass bracket is defined to be $\{f,g\}_{m}=\Im(f\overline{g})$. Note that in the particular case
when $\mathcal{N}=F'(|\phi|^{2})\phi$ one has $M_{a}^0=\partial_{t}V_{a}$.

Another calculation establishes
\begin{lemma}
$$\partial_{t}M_{a}^0=\int_{\R^n} (-\Delta\Delta a)|\phi|^{2}+4\int_{\R^n}
a_{jk}\Re(\overline{\phi_{j}}\phi_{k})+2\int_{\R^n}
a_{j}\{\mathcal{N},\phi\}_{p}^{j},$$ where we define the momentum
bracket to be
$\{f,g\}_p=\Re(f\nabla\overline{g}-g\nabla\overline{f})$ and repeated indices are implicitly summed.
\end{lemma}

Note again that when $\mathcal{N}=F'(|\phi|^{2})\phi$ we have
$\{\mathcal{N},\phi\}_p=-\nabla G(|\phi|^{2})$ for $G(x)=xF'(x)-F(x)$.
In particular, in the cubic case
$\{\mathcal{N},\phi\}_p=-\frac{1}{2}\nabla(|\phi|^{4})$.

Now let $a(x)=|x|$. Easy computations show that for dimension
$n\geq 4$ we have the following identities:
\begin{align*}
a_{j}(x)=&\frac{x_{j}}{|x|} \\
a_{jk}(x)=&\frac{\delta_{jk}}{|x|}-\frac{x_{j}x_{k}}{|x|^{3}} \\
\Delta a(x)=&\frac{n-1}{|x|} \\
-\Delta \Delta a(x)=&\frac{(n-1)(n-3)}{|x|^{3}}.
\end{align*}

For this choice of the function $a$, one should interpret the $M_a^0$ as a spatial average of
the radial component of the $L^2_x$-mass current. Taking its time derivative we get
\begin{align*}
\partial_{t}M_{a}^0
&=(n-1)(n-3)\int_{\R^n} \frac{|\phi(x)|^{2}}{|x|^{3}}dx
   +4\int_{\R^n} (\frac{\delta_{jk}}{|x|}-\frac{x_{j}x_{k}}{|x|^{3}}) \Re(\overline{\phi_{j}}\phi_{k})(x)dx\\
&\quad +2\int_{\R^n} \frac{x_{j}}{|x|}\{\mathcal{N},\phi\}_{p}^{j}(x)dx \\
&=(n-1)(n-3)\int_{\R^n} \frac{|\phi(x)|^{2}}{|x|^{3}}dx
   +4\int_{\R^n} \frac{1}{|x|} |\nabla_{0}\phi(x)|^{2}dx\\
&\quad +2\int_{\R^n} \frac{x}{|x|} \{\mathcal{N},\phi\}_{p}(x)dx,
\end{align*}
where we use $\nabla_{0}$ to denote the complement of the radial
portion of the gradient, that is
$\nabla_{0}=\nabla-\frac{x}{|x|}(\frac{x}{|x|}\nabla)$.

We may center the above argument at any other point $y\in
\R^{4}$. Choosing $a(x)=|x-y|$, we define the Morawetz action centered at $y$ to be
$$
M_{a}^y(t)=2\int_{\R^{n}}\frac{x-y}{|x-y|}\Im(\overline{\phi(x)}\nabla \phi(x))dx.
$$
The same computations now yield
\begin{align*}
\partial_{t}M_a^{y}
&=(n-1)(n-3)\int_{\R^n} \frac{|\phi(x)|^{2}}{|x-y|^{3}}dx
   +4\int_{\R^n} \frac{1}{|x-y|} |\nabla_{y}\phi(x)|^{2}dx\\
&\quad+2\int_{\R^n} \frac{x-y}{|x-y|} \{\mathcal{N},\phi\}_{p}(x)dx.
\end{align*}

We are now ready to define the interaction Morawetz potential, which is a way of quantifying
how mass is interacting with (moving away from) itself:
\begin{align*}
M^{interact}(t)
&=\int_{\R^n} |\phi(t,y)|^{2}M_a^{y}(t)dy\\
&=2\Im \int_{\R^n} \int_{\R^n} |\phi(t,y)|^{2}\frac{x-y}{|x-y|}\nabla\phi(t,x)\overline{\phi(t,x)}dxdy.
\end{align*}
One gets immediately the easy estimate
$$|M^{interact}(t)|
   \leq 2 \|\phi(t)\|_{L_{x}^{2}}^{3} \|\phi(t)\|_{\dot{H}_{x}^{1}}.$$

Calculating the time derivative of the interaction Morawetz
potential we get the following virial-type identity,
\begin{align}
\partial_{t}M^{interact}
=&(n-1)(n-3)\int_{\R^n} \int_{\R^n} \frac{|\phi(y)|^{2}|\phi(y)|^{2}}{|x-y|^{3}}dxdy  \label{vi1} \\
 &+4\int_{\R^n} \int_{\R^n} \frac{|\phi(y)|^{2}|\nabla_{y}\phi(x)|^{2}}{|x-y|}dxdy  \label{vi2} \\
 &+2\int_{\R^n} \int_{\R^n} \frac{|\phi(y)|^{2}}{|x-y|}(x-y) \{\mathcal{N},\phi\}_{p}(x)dxdy  \label{vi3} \\
 &+2 \int_{\R^n} \partial_{y_{k}} \Im(\phi\overline{\phi_{k}})(y)M_a^{y}dy  \label{vi4} \\
 &+4\Im \int_{\R^n} \int_{\R^n} \{\mathcal{N},\phi\}_{m}(y)\frac{x-y}{|x-y|}\nabla\phi(x)\overline{\phi(x)}dxdy. \label{vi5}
\end{align}

As far as the terms in the above virial-type identity are concerned, we will establish
\begin{lemma}\label{lemmaviterms}
\eqref{vi4} $\geq$ --\eqref{vi2}.
\end{lemma}

Thus, integrating over the compact interval $I_{0}$ we get

\begin{proposition}[Interaction Morawetz inequality] \label{intmorineq}
\begin{align*}
(n-1)(n-3)&\int_{I_{0}} \int_{\R^n} \int_{\R^n} \frac{|\phi(t,y)|^{2}|\phi(t,x)|^{2}}{|x-y|^{3}}dxdydt\\
 &+2\int_{I_{0}} \int_{\R^n} \int_{\R^n} \frac{|\phi(t,y)|^{2}}{|x-y|}(x-y)\{\mathcal{N},\phi\}_{p}(t,x)dxdydt \\
  \leq {}&2\|\phi\|_{L_{t}^{\infty}L_{x}^{2}(I_{0}\times\R^{n})}^{3} \|\phi\|_{L_{t}^{\infty}\dot{H}_{x}^{1}(I_{0}\times\R^{n})} \\
          &+4\int_{I_{0}}\int_{\R^n}\int_{\R^n}|\{\mathcal{N},\phi\}_{m}(t,y)||\nabla\phi(t,x)||\phi(t,x)|dxdydt.
\end{align*}
\end{proposition}

Note that in the particular case $\mathcal{N}=|u|^{2}u$, after performing an integration by parts in
the momentum bracket term, the inequality becomes
\begin{align*}
(n-1)(n-3)&\int_{I_{0}} \int_{\R^n} \int_{\R^n} \frac{|u(t,y)|^{2}|u(t,x)|^{2}}{|x-y|^{3}}dxdydt\\
 &\qquad+(n-1)\int_{I_{0}} \int_{\R^n} \int_{\R^n} \frac{|u(t,y)|^{2}|u(t,x)|^{4}}{|x-y|}dxdydt\\
&\leq 2\|u\|_{L_{t}^{\infty}L_{x}^{2}(I_{0}\times\R^{n})}^{3} \|u\|_{L_{t}^{\infty}\dot{H}_{x}^{1}(I_{0}\times\R^{n})}.
\end{align*}

We turn now to the proof of Lemma \ref{lemmaviterms}.
We write
$$
\eqref{vi4}=4\int_{\R^n} \int_{\R^n} \partial_{y_{k}} \Im(\phi(y)\overline{\phi_{k}(y)})\frac{x_j-y_j}{|x-y|}\Im(\overline{\phi(x)}\phi_j(x))dxdy,
$$
where we sum over repeated indices. We integrate by parts moving $\partial_{y_k}$ to the unit vector
$\frac{x-y}{|x-y|}$. Using the identity
$$
\partial_{y_k}\bigl(\frac{x_j-y_j}{|x-y|}\bigr)=-\frac{\delta_{kj}}{|x-y|}+\frac{(x_k-y_k)(x_j-y_j)}{|x-y|^3}
$$
and the notation $p(x)=2\Im(\overline{\phi(x)}\nabla \phi(x))$ for the momentum density, we rewrite \eqref{vi4} as
$$
-\int_{\R^n} \int_{\R^n} \Bigl[ p(y)p(x)-\Bigl( p(y)\frac{x-y}{|x-y|} \Bigr) \Bigl( p(x)\frac{x-y}{|x-y|} \Bigr) \Bigr] \frac{dxdy}{|x-y|}.
$$
Note that the quantity between the square brackets represents the inner product between the projection of the momentum density
$p(y)$ onto the orthogonal complement of $(x-y)$ and the projection of $p(x)$ onto the same space. But
\begin{align*}
|\pi_{(x-y)^\perp}p(y)|
   &=|p(y)-\frac{x-y}{|x-y|}(\frac{x-y}{|x-y|}p(y))|
    =2|\Im(\overline{\phi(y)}\nabla_x \phi(y))| \\
   &\leq 2|\phi(y)||\nabla_x \phi(y))|.
\end{align*}
As the same estimate holds when we switch $y$ and $x$, we get
\begin{align*}
\eqref{vi4}
    &\geq -4\int_{\R^n} \int_{\R^n} |\phi(y)||\nabla_x \phi(y))| |\phi(x)| |\nabla_y \phi(x))|dxdy \\
    &\geq -2 \int_{\R^n} \int_{\R^n} \frac{|\phi(y)|^{2}|\nabla_{y}\phi(x)|^{2}}{|x-y|}dxdy
       -2 \int_{\R^n} \int_{\R^n} \frac{|\phi(x)|^{2}|\nabla_{x}\phi(y)|^{2}}{|x-y|}dxdy \\
    &\geq -\eqref{vi2}.
\end{align*}

\subsection{Morawetz inequality: the setup}
We are now ready to start the proof of Proposition \ref{prop flim}. As the
statement is invariant under scaling, we normalize $N_{*}=1$
and define $u_{hi}=P_{> 1}u$ and $u_{lo}=P_{\leq 1}u$. As
we assume $1=N_{*}<c(\eta_{1})N_{min}$, we get
$1<c(\eta_{1})N(t)$, $\forall t\in I_{0}$. Provided we choose
$c(\eta_{1})$ sufficiently small, the frequency localization
result and Sobolev yield
\begin{equation}\label{sijz}
\|u_{<\eta_1^{-1}}\|_{L_{t}^{\infty}\dot{H}_{x}^{1}(I_{0}\times\R^{4})}
  +\|u_{<\eta_1^{-1}}\|_{L_{t}^{\infty}L_{x}^{4}(I_{0}\times\R^{4})}
  \leq \eta_{1}.
\end{equation}
Hence, $u_{lo}$ has small energy
\begin{align}
\|u_{lo}\|_{L_{t}^{\infty}\dot{H}_{x}^{1}(I_{0}\times\R^{4})}
  +\|u_{lo}\|_{L_{t}^{\infty}L_{x}^{4}(I_{0}\times\R^{4})}\lesssim
\eta_{1}. \label{lowfreqsmall}
\end{align}
Using \eqref{mass high freq bound} and \eqref{sijz}, one also sees that $u_{hi}$ has small mass
\begin{align}
\|u_{hi}\|_{L_{t}^{\infty}L_{x}^{2}(I_{0}\times\R^{4})}\lesssim \eta_{1}. \label{smallmass}
\end{align}

Our goal is to prove
\begin{align}\label{pfu}
\int_{I_{0}}\int_{\R^{4}}\int_{\R^{4}}\frac{|u_{hi}(t,x)|^{2}|u_{hi}(t,y)|^{2}}{|x-y|^{3}}dx dy dt\lesssim\eta_{1}.
\end{align}

Since in four dimensions convolution with $1/|x|^3$ is basically the same as the
fractional integration operator $|\nabla|^{-1}$, the above estimate translates into
\begin{align}
\| |u_{hi}|^2 \|_{L^2_t \dot H^{-1/2}_x(I_{0} \times \R^4)} \lesssim \eta^{1/2}_1. \label{hifreqasmp}
\end{align}

By a standard continuity argument, it will suffice to prove \eqref{hifreqasmp} under the bootstrap hypothesis
\begin{align}
\| |u_{hi}|^2 \|_{L^2_t \dot H^{-1/2}_x(I_{0} \times \R^4)} \leq C_0^{\frac{1}{2}} \eta^{1/2}_1, \label{bootstraphyp1}
\end{align}
for a large constant $C_0$ depending on energy but not on any of the $\eta$'s. More rigorously, one needs to prove that
\eqref{bootstraphyp1} implies \eqref{pfu} whenever $I_0$ is replaced by a subinterval of $I_0$ in order to run
the bootstrap argument correctly. However, it will become clear to the reader that the argument below works not only
for $I_0$ but also for any subinterval of $I_0$.

Note that a consequence of \eqref{bootstraphyp1} is
\begin{align}
\|P_{\leq N}|u_{hi}|^2 \|_{L^2_{t,x}(I_0\times \R^4)}\leq N^{\frac{1}{2}} C_0^{\frac{1}{2}}\eta^{1/2}_1. \label{bootstraphyp2}
\end{align}

\begin{proposition} \label{etprop}
With the notation and assumptions above we have
$$\int_{I_{0}}\int_{\R^4}\int_{\R^4}
\frac{|u_{hi}(t,y)|^{2}|u_{hi}(t,x)|^{2}}{|x-y|^{3}}dxdydt+\int_{I_{0}}\int_{\R^4}\int_{\R^4}
\frac{|u_{hi}(t,y)|^{2}|u_{hi}(t,x)|^{4}}{|x-y|}dxdydt
$$
\begin{align}
\lesssim& \eta_1^3 \label{et0}\\
&+\eta_{1} \sum_{j=0}^2 \int_{I_{0}} \int_{\R^4}
    |P_{hi}\O(u_{hi}^{j}u_{lo}^{3-j})(t,y)||u_{hi}(t,y)|dydt \label{et1}\\
&+\eta_1 \int_{I_{0}}\int_{\R^4}|P_{lo}\O(u_{hi}^{3})(t,y)||u_{hi}(t,y)|dydt \label{et2}\\
&+\eta_{1}^{2}\sum_{j=1}^{3}\int_{I_{0}}\int_{\R^4}|u_{hi}(t,x)|^{j}|u_{lo}(t,x)|^{3-j}|\nabla u_{lo}(t,x)|dxdt \label{et3}\\
&+\eta_{1}^{2}\int_{I_{0}}\int_{\R^4}|u_{hi}\nabla P_{lo}(|u|^{2}u)|(t,x)dxdt \label{et4}\\
&+\sum_{j=1}^{3}\Bigl| \int_{I_{0}}\int_{\R^4}\int_{\R^4}\frac{|u_{hi}(t,y)|^{2}\O(u_{hi}^{j}u_{lo}^{4-j})(t,x)}{|x-y|}dxdydt\Bigr|. \label{et5}
\end{align}

\end{proposition}

\begin{proof}

Applying Proposition \ref{intmorineq} with $\phi=u_{hi}$ and
$\mathcal{N}=P_{hi}(|u|^{2}u)$ we find
\begin{align*}
3\int_{I_{0}} \int_{\R^4} \int_{\R^4} & \frac{|u_{hi}(t,y)|^{2}|u_{hi}(t,x)|^{2}}{|x-y|^{3}}dxdydt\\
&\quad+2\int_{I_{0}} \int_{\R^4} \int_{\R^4} \frac{|u_{hi}(t,y)|^{2}}{|x-y|}(x-y)\{P_{hi}(|u|^{2}u),u_{hi}\}_{p}(t,x)dxdydt \\
   &\leq 2\|u_{hi}\|_{L_{t}^{\infty}L_{x}^{2}(I_{0}\times\R^{4})}^{3} \|u_{hi}\|_{L_{t}^{\infty}\dot{H}_{x}^{1}(I_{0}\times\R^{4})} \\
   &\quad +4 \int_{I_{0}} \int_{\R^4} \int_{\R^4} |\{P_{hi}(|u|^{2}u),u_{hi}\}_{m}(t,y)| |\nabla u_{hi}(t,x)||u_{hi}(t,x)|dxdydt.
\end{align*}

Observe that \eqref{smallmass} plus conservation of energy dictates
$$
\|u_{hi}\|_{L_{t}^{\infty}L_{x}^{2}(I_{0}\times\R^{4})}^{3}\|u_{hi}\|_{L_{t}^{\infty}\dot{H}_{x}^{1}(I_{0}\times\R^{4})}
   \lesssim \eta^3_1,
$$
which is the error term \eqref{et0}.\\

We deal first with the mass bracket term. Exploiting cancellation we write
$$\{P_{hi}(|u|^{2}u),u_{hi}\}_{m}=\{P_{hi}(|u|^{2}u)-|u_{hi}|^{2}u_{hi},u_{hi}\}_{m}.$$

Writing
\begin{align*}
P_{hi}(|u|^{2}u)-|u_{hi}|^{2}u_{hi}
  &=P_{hi}(|u|^{2}u-|u_{hi}|^{2}u_{hi})-P_{lo}(|u_{hi}|^{2}u_{hi}) \\
  &=\sum_{j=0}^{2}P_{hi}\O(u_{hi}^{j}u_{lo}^{3-j})-P_{lo}\O(u_{hi}^{3}),
\end{align*}
we estimate the mass bracket term by the error terms \eqref{et1} and \eqref{et2} as follows:
\begin{align*}
\sum_{j=0}^2 \int_{I_{0}} \int_{\R^4} \int_{\R^4} & |u_{hi}(t,x)||\nabla u_{hi}(t,x)||P_{hi}\O(u_{hi}^{j}u_{lo}^{3-j})(t,y)||u_{hi}(t,y)|dxdydt \\
  &+\int_{I_{0}} \int_{\R^4} \int_{\R^4}|u_{hi}(t,x)||\nabla u_{hi}(t,x)||P_{lo}\O(u_{hi}^{3})(t,y)||u_{hi}(t,y)|dxdydt \\
    \lesssim &\eta_{1}\sum_{j=0}^2 \int_{I_{0}} \int_{\R^4} |P_{hi}\O(u_{hi}^{j}u_{lo}^{3-j})(t,y)||u_{hi}(t,y)|dydt \\
             &+\eta_{1}\int_{I_{0}}\int_{\R^4} |P_{lo}\O(u_{hi}^{3})(t,y)||u_{hi}(t,y)|dydt,
\end{align*}
where in order to obtain the last inequality we used
$$
\int |u_{hi}(t,x)||\nabla u_{hi}(t,x)|dx\lesssim \|u_{hi}\|_{L_{t}^{\infty}L_{x}^{2}(I_{0}\times\R^{4})}
\|\nabla u_{hi}\|_{L_{t}^{\infty}L_{x}^{2}(I_{0}\times\R^{4})}
   \lesssim \eta_1.
$$

We turn now towards the momentum bracket term and write
\begin{align*}
\{P_{hi}(|u|^{2}u),u_{hi}\}_{p}
 &=\{|u|^{2}u,u_{hi}\}_{p}-\{P_{lo}(|u|^{2}u),u_{hi}\}_{p} \\
 &=\{|u|^{2}u,u\}_{p}-\{|u|^{2}u,u_{lo}\}_{p}-\{P_{lo}(|u|^{2}u),u_{hi}\}_{p} \\
 &=\{|u|^{2}u,u\}_{p}-\{|u_{lo}|^{2}u_{lo},u_{lo}\}_{p}-\{|u|^{2}u-|u_{lo}|^{2}u_{lo},u_{lo}\}_{p}\\
 &\quad -\{P_{lo}(|u|^{2}u),u_{hi}\}_{p} \\
 &=-\frac{1}{2}\nabla(|u|^{4}-|u_{lo}|^{4})-\{|u|^{2}u-|u_{lo}|^{2}u_{lo},u_{lo}\}_{p}\\
 &\quad-\{P_{lo}(|u|^{2}u),u_{hi}\}_{p}.
\end{align*}

After an integration by parts, the first term in the momentum
bracket contributes a multiple of
\begin{align*}
\int_{I_{0}} \int_{\R^4} \int_{\R^4} &\frac{|u_{hi}(t,y)|^{2}|u_{hi}(t,x)|^{4}}{|x-y|}dxdydt \\
                        &+\sum_{j=1}^{3}\Bigl|\int_{I_{0}} \int_{\R^4} \int_{\R^4} \frac{|u_{hi}(t,y)|^{2}\O(u_{hi}^{j}u_{lo}^{4-j})(t,x)}{|x-y|}dxdydt\Bigr|,
\end{align*}
in which we recognize \eqref{et5} and the left-hand side term in Proposition \ref{etprop}.\\

In order to estimate the contribution of the second term in the
momentum bracket, we write $\{f,g\}_{p}=\nabla \O(fg)+\O(f\nabla g)$ and
$$|u|^{2}u-|u_{lo}|^{2}u_{lo}=\sum_{j=1}^{3}\O(u_{hi}^{j}u_{lo}^{3-j})$$ and we decompose
\begin{equation*}
\{|u|^{2}u-|u_{lo}|^{2}u_{lo},u_{lo}\}_{p}
   =\sum_{j=1}^{3} \nabla\O(u_{hi}^{j}u_{lo}^{4-j})
    +\sum_{j=1}^{3} \O(u_{hi}^{j}u_{lo}^{3-j}\nabla u_{lo}) = I + II.
\end{equation*}

In order to estimate the contribution of $I$ to the momentum bracket term, we
integrate by parts to rediscover \eqref{et5}.

Taking absolute values inside the integrals, $II$ contributes to the momentum bracket term
\begin{align*}
\sum_{j=1}^{3} \int_{I_{0}} &\int_{\R^4} \int_{\R^4} |u_{hi}(t,y)|^{2}|u_{hi}(t,x)|^{j}|u_{lo}(t,x)|^{3-j}|\nabla u_{lo}(t,x)|dxdydt \\
                                      &\lesssim \|u_{hi}\|_{L_{t}^{\infty}L_{x}^{2}(I_{0}\times\R^{4})}^{2} \sum_{j=1}^{3}\int_{I_{0}} \int_{\R^4} |u_{hi}(t,x)|^{j}|u_{lo}(t,x)|^{3-j}|\nabla u_{lo}(t,x)|dxdt \\
                                      &\lesssim \eta_{1}^{2} \sum_{j=1}^{3} \int_{I_{0}} \int_{\R^4} |u_{hi}(t,x)|^{j}|u_{lo}(t,x)|^{3-j}|\nabla u_{lo}(t,x)|dxdt,
\end{align*} which is \eqref{et3}.

The only term left to consider is the last term in the momentum
bracket. When the derivative (from the definition of the momentum
bracket) falls on $P_{lo}(|u|^{2}u)$, we estimate its final
contribution by taking absolute values inside the integrals to get
\begin{align*}
\int_{I_{0}} \int_{\R^4} \int_{\R^4} &|u_{hi}(t,y)|^{2}|u_{hi}\nabla P_{lo}(|u|^{2}u)|(t,x)dxdydt \\
                       &\lesssim \|u_{hi}\|_{L_{t}^{\infty}L_{x}^{2}(I_{0}\times\R^{4})}^{2} \int_{I_{0}} \int_{\R^4} |u_{hi}\nabla P_{lo}(|u|^{2}u)|(t,x)dxdt \\
                       &\lesssim \eta_{1}^{2} \int_{I_{0}} \int_{\R^4} |u_{hi}\nabla P_{lo}(|u|^{2}u)|(t,x)dxdt,
\end{align*}
in which we recognize the error term \eqref{et4}.

When the derivative falls on $u_{hi}$, we first integrate by parts
to get as a final contribution
\begin{align}
\int_{I_{0}} \int_{\R^4} \int_{\R^4} &|u_{hi}(t,y)|^{2}|u_{hi}\nabla P_{lo}(|u|^{2}u)|(t,x)dxdydt \\
                       &+\Bigl|\int_{I_{0}} \int_{\R^4} \int_{\R^4} \frac{|u_{hi}(t,y)|^{2}u_{hi}(t,x)P_{lo}(|u|^{2}u)(t,x)}{|x-y|}dxdydt\Bigr| \\
              \lesssim {}&\eta_{1}^{2}\int_{I_{0}} \int_{\R^4} |u_{hi}\nabla P_{lo}(|u|^{2}u)|(t,x)dxdt \\
                       &+\Bigl|\int_{I_{0}} \int_{\R^4} \int_{\R^4} \frac{|u_{hi}(t,y)|^{2}u_{hi}(t,x)P_{lo}(|u|^{2}u)(t,x)}{|x-y|}dxdydt\Bigr|. \label{a3}
\end{align}

Decomposing $u=u_{lo}+u_{hi}$, we bound \eqref{a3} by \eqref{et5} and the second term on the left-hand side of
Proposition \ref{etprop}.

\end{proof}

\subsection{Strichartz control on low and high frequencies} The
purpose of this section is to obtain estimates on the low and high
frequency parts of $u$, which we will use to bound the error
terms obtained in the previous section.

\begin{proposition}[Strichartz control on low frequencies] \label{Slf}
We have
\begin{align}
\|u_{lo}\|_{\dot{S}^{1}(I_0\times \R^4)}\lesssim \eta_{1}^{\frac{1}{2}}C_{0}^{\frac{1}{2}}. \label{lowfreqcontrol}
\end{align}
\end{proposition}

\begin{proof}
Throughout the proof all spacetime norms will be on $I_0\times \R^4$.  Lemma \ref{lemma linear strichartz} dictates
$$
\|u_{lo}\|_{\dot{S}_{1}} \lesssim \|u_{lo}(t_0)\|_{\dot{H}_{x}^1}
     +\|\nabla P_{lo}(|u|^{2}u)\|_{L_{t}^2L_x^{4/3}}.
$$
Note that according to our assumptions,
$\|u_{lo}(t_0)\|_{\dot{H}_{x}^1}\lesssim \eta_1$.

Write $\nabla P_{lo}(|u|^{2}u)=\sum_{j=0}^{3}\nabla P_{lo}\O(u_{hi}^{j}u_{lo}^{3-j})$. The contribution coming from
the $j=0$ term can be estimated by
\begin{align*}
\|\nabla P_{lo}\O(u_{lo}^3)\|_{L_{t}^2L_x^{4/3}} \lesssim \|u_{lo}^2 \nabla u_{lo}\|_{L_{t}^2L_x^{4/3}}
    \lesssim \|\nabla u_{lo}\|_{L^2_t L^4_x} \|u_{lo}\|_{L^{\infty}_tL^4_x}^{2}
    \lesssim \eta_1^2 \|u_{lo}\|_{\dot{S}_{1}}.
\end{align*}
We estimate the contribution coming from the term corresponding to $j=1$ using Bernstein and \eqref{smallmass} as
$$
\|\nabla P_{lo}\O(u_{lo}^2u_{hi})\|_{L_{t}^2L_x^{4/3}}
   \lesssim \|u_{lo}^2u_{hi}\|_{L_{t}^2L_x^{4/3}}
   \lesssim\|u_{lo}\|_{L^4_tL^8_x}^{2} \|u_{hi}\|_{L^{\infty}_t L^2_x}
   \lesssim \eta_{1} \|u_{lo}\|_{\dot{S}_{1}}^{2}.
$$
The term corresponding to $j=2$ is estimated via Bernstein and \eqref{bootstraphyp2} as
\begin{align*}
\|\nabla P_{lo}\O(u_{lo} u_{hi}^{2})\|_{L_t^2L_x^{4/3}}
&\lesssim\|P_{lo}\O(u_{lo} u_{hi}^{2})\|_{L_t^2L_x^{4/3}}
   \lesssim \|u_{lo}\|_{L^{\infty}_t L^4_x} \|P_{\leq 10}(u_{hi}^{2})\|_{L^2_{t,x}}\\
&\lesssim \eta_1^{\frac{3}{2}} C_{0}^{\frac{1}{2}}.
\end{align*}

We now turn toward the $j=3$ term; consider a dyadic piece
$\nabla P_{lo}\O(u_{N_{1}} u_{N_{2}}u_{N_{3}})$ for $N_1\geq N_2\geq N_3\geq 1$, and notice that we can replace
$u_{N_{2}} u_{N_{3}}$ by $P_{\leq 10 N_1}(u_{N_{2}}u_{N_{3}})$. Using Bernstein and \eqref{bootstraphyp2} we estimate
\begin{align*}
\|\nabla P_{lo}\O(u_{hi}^{3})\|_{L_t^2L_x^{4/3}}
       &\lesssim\|P_{lo}\O(u_{hi}^{3})\|_{L_t^2L_x^1}\\
       &\lesssim \sum_{N_1\geq N_2\geq N_3\geq 1} \|u_{N_{1}}\|_{L^{\infty}_t L^2_x} \|P_{\leq 10 N_1}(u_{N_{2}} u_{N_{3}})\|_{L^2_{t,x}} \\
       &\lesssim \sum_{N_1\geq N_2\geq N_3\geq 1} N_{1}^{-1}N_{1}^{\frac{1}{2}} \eta_{1}^{\frac{1}{2}} C_{0}^{\frac{1}{2}}\\
       &\lesssim \eta_{1}^{\frac{1}{2}} C_{0}^{\frac{1}{2}}.
\end{align*}

Putting everything together we obtain
$$
\|u_{lo}\|_{\dot{S}_{1}}
  \lesssim \eta_1 + \eta_1^2 \|u_{lo}\|_{\dot{S}_{1}} + \eta_1 \|u_{lo}\|_{\dot{S}_{1}}
           + \eta_1^{\frac{3}{2}} C_{0}^{\frac{1}{2}} + \eta_1^{\frac{1}{2}} C_{0}^{\frac{1}{2}}.
$$
Choosing $\eta_1$ sufficiently small, we get $\|u_{lo}\|_{\dot{S}_{1}}\lesssim \eta_1^{\frac{1}{2}} C_{0}^{\frac{1}{2}}.$
\end{proof}

\begin{proposition}[Strichartz control on high frequencies] We
have
\begin{align}
\|\nabla^{-\frac{1}{2}}u_{hi}\|_{L^2_t L^4_x(I_0\times\R^4)} \lesssim C_{0}^{\frac{1}{2}} \eta_1^{\frac{1}{2}}. \label{highfreqcontrol}
\end{align}
\end{proposition}

\begin{proof} Throughout the proof all spacetime norms will be on $I_0\times\R^4$. Strichartz dictates
$$
\|\nabla ^{-\frac{1}{2}}u_{hi}\|_{L^2_tL^4_x}
   \lesssim \|\nabla^{-\frac{1}{2}}u_{hi}\|_{L^{\infty}_tL^2_x} + \|\nabla^{-\frac{1}{2}}P_{hi}(|u|^{2}u)\|_{L^2_tL^{4/3}_x}.
$$

We estimate the first term crudely by
$$\|\nabla^{-\frac{1}{2}}u_{hi}\|_{L^{\infty}_tL^2_x} \lesssim \|u_{hi}\|_{L^{\infty}_tL^2_x}\lesssim \eta_1.$$

To deal with the nonlinearity we decompose $u$ into $u_{lo}$ and $u_{hi}$, and write
$\nabla^{-\frac{1}{2}} P_{hi}(|u|^{2}u)=\sum_{j=0}^{3}\nabla^{-\frac{1}{2}} P_{hi}\O(u_{hi}^{j}u_{lo}^{3-j})$.
The control we have gained on low frequencies in Proposition~\ref{Slf} and Bernstein allow us estimate the terms
corresponding to the cases $j=0$ and $j=1$:
\begin{align*}
\|\nabla^{-\frac{1}{2}}P_{hi}\O(u_{lo}^3)\|_{L^2_tL^{4/3}_x}
&\lesssim \|\nabla P_{hi}\O(u_{lo}^3)\|_{L^2_tL^{4/3}_x}
 \lesssim  \|\nabla u_{lo}\|_{L^2_tL^4_x} \|u_{lo}\|^2_{L^{\infty}_tL^4_x}\\
&\lesssim \eta_1^2 \|u_{lo}\|_{\dot{S}^1}
 \lesssim C_{0}^{\frac{1}{2}}\eta_1^{\frac{5}{2}} \\
\|\nabla^{-\frac{1}{2}}P_{hi}\O(u_{hi} u_{lo}^2)\|_{L^2_tL^{4/3}_x}
&\lesssim \|P_{hi}\O(u_{hi} u_{lo}^2)\|_{L^2_tL^{4/3}_x}
 \lesssim \|u_{hi}\|_{L_t^\infty L_x^2}\|u_{lo}\|_{L_t^4L_x^8}^2\\
&\lesssim \eta_1 \|u_{lo}\|_{\dot{S}^1}^2
 \lesssim C_0\eta_1^2.
\end{align*}

In order to deal with the terms corresponding to the cases $j=2$ and $j=3$, we will prove the more general
estimate
\begin{align}
\|\nabla^{-\frac{1}{2}}\{(\nabla^{\frac{1}{2}} f)(\nabla^{-\frac{1}{2}}g)\}\|_{L_x^{4/3}}
   \lesssim \|f\|_{L_{x}^{2}}\|g\|_{L_x^{8/3}}. \label{paraproduct}
\end{align}

In order to prove \eqref{paraproduct}, we decompose the left-hand side into $\pi_{h,h}$, $\pi_{l,h}$, and $\pi_{h,l}$
which represent the projections onto high-high, low-high, and high-low frequency interactions.  More precisely, for any pair of functions $(\phi, \psi)$,
we write $\pi_{h,h}(\phi, \psi)=\sum_{N\sim M}P_N \phi P_M \psi$, $\pi_{l,h}(\phi, \psi)=\sum_{N\lesssim M}P_N \phi P_M \psi$,
and $\pi_{h,l}(\phi, \psi)=\sum_{N\gtrsim M}P_N \phi P_M \psi$.

The high-high and low-high frequency interactions are going to be treated in the same manner. Let's consider for
example the first one. A simple application of Sobolev embedding yields
\begin{align}
\|\nabla^{-\frac{1}{2}}\pi_{h,h}\{(\nabla^{\frac{1}{2}} f)(\nabla^{-\frac{1}{2}}g)\}\|_{L_x^{4/3}}
   \lesssim \|\pi_{h,h}\{(\nabla^{\frac{1}{2}}f)(\nabla^{-\frac{1}{2}} g)\}\|_{L_x^{8/7}}. \label{goodlemma}
\end{align}
Now we only have to notice that the multiplier associated to the operator
$T(f,g)=\pi_{h,h}\{(\nabla^{\frac{1}{2}} f)(\nabla^{-\frac{1}{2}}g)\}$, i.e.
$$
\sum_{N\sim M} |\xi_1|^{\frac{1}{2}}\widehat{P_N f}(\xi_1) |\xi_2|^{-\frac{1}{2}}\widehat{P_M g}(\xi_2),
$$
is a symbol of order one with $\xi = (\xi_1,\xi_2)$,
since then a theorem of Coifman and Meyer (\cite{coifmey:1}, \cite{coifmey:2}) will conclude our claim.

To deal with the $\pi_{h,l}$ term, we first notice that the multiplier associated to the operator
$T(f,h)=\nabla^{-\frac{1}{2}} \pi_{h,l}\{(\nabla^{\frac{1}{2}} f)h\}$, i.e.
$$
\sum_{N\gtrsim M}|\xi_1+\xi_2|^{-\frac{1}{2}}|\xi_1|^{\frac{1}{2}}\widehat{P_N f}(\xi_1) \widehat{P_M h}(\xi_2),
$$
is an order one symbol. The result cited above yields
$$
\|\nabla^{-\frac{1}{2}}\pi_{h,l}\{(\nabla^{\frac{1}{2}} f)(\nabla^{-\frac{1}{2}}g)\}\|_{L_x^{4/3}}
     \lesssim \|f\|_{L^2_x} \|\nabla^{-\frac{1}{2}}g\|_{L_x^4}.
$$
Finally, Sobolev embedding dictates the estimate $\|\nabla^{-\frac{1}{2}}g\|_{L_x^4}\lesssim \|g\|_{L_x^{8/3}}$.

To estimate the contribution of the term corresponding to the case $j=2$, we take $f=\nabla^{-\frac{1}{2}}|u_{hi}|^2$
and $g=\nabla^{\frac{1}{2}}u_{lo}$ in \eqref{paraproduct}. Using Sobolev embedding and \eqref{bootstraphyp1} we get
\begin{align*}
\|\nabla^{-\frac{1}{2}}P_{hi}\O(u_{hi}^2u_{lo})\|_{L^2_tL^{4/3}_x}
&\lesssim \|\nabla^{-\frac{1}{2}}|u_{hi}|^2\|_{L^2_{t,x}} \|\nabla^{\frac{1}{2}}u_{lo}\|_{L^{\infty}_tL_x^{8/3}}
     \lesssim \eta_1^{\frac{1}{2}} C_0^{\frac{1}{2}} \|\nabla u_{lo}\|_{L^{\infty}_tL_x^2}\\
&\lesssim C_0^{\frac{1}{2}}\eta_1^{\frac{1}{2}}\eta_1.
\end{align*}
Similarly, to treat the case $j=3$ we take $f=\nabla^{-\frac{1}{2}}|u_{hi}|^2$ and $g=\nabla^{\frac{1}{2}}u_{hi}$
in \eqref{paraproduct} and use Sobolev embedding and \eqref{bootstraphyp1} to estimate
\begin{align*}
\|\nabla^{-\frac{1}{2}}P_{hi}\O(u_{hi}^3)\|_{L^2_tL^{4/3}_x}
&\lesssim \|\nabla^{-\frac{1}{2}}|u_{hi}|^2\|_{L^2_{t,x}} \|\nabla^{\frac{1}{2}}u_{hi}\|_{L^{\infty}_tL_x^{8/3}}
     \lesssim \eta_1^{\frac{1}{2}} C_0^{\frac{1}{2}} \|\nabla u_{hi}\|_{L^{\infty}_tL_x^2}\\
&\lesssim \eta_1^{\frac{1}{2}} C_0^{\frac{1}{2}}.
\end{align*}

Combining the above estimates we get
$$
\|\nabla ^{-\frac{1}{2}}u_{hi}\|_{L^2_tL^4_x}
    \lesssim \eta_1 + C_{0}^{\frac{1}{2}}\eta_1^{\frac{5}{2}} +C_0\eta_1^2 + C_0^{\frac{1}{2}}\eta_1^{\frac{3}{2}} + \eta_1^{\frac{1}{2}} C_0^{\frac{1}{2}}
    \lesssim \eta_1^{\frac{1}{2}} C_0^{\frac{1}{2}}.
$$
which proves the proposition.
\end{proof}

\begin{corollary}\label{cor uhi}
\begin{align}
\|u_{hi}\|_{L^3_{t,x}(I_0\times\R^4)} \lesssim C_{0}^{\frac{1}{3}} \eta_1^{\frac{1}{3}}. \label{u_hi estimate}
\end{align}
\end{corollary}

\begin{proof} The claim follows interpolating between $$\nabla u_{hi}\in L^{\infty}_t L^2_x(I_0\times\R^4)$$ and
$$\nabla^{-\frac{1}{2}} u_{hi}\in L^2_tL^4_x(I_0\times\R^4).$$
\end{proof}

\subsection{Morawetz inequality: error terms} In this section we use the control on $u_{lo}$ and $u_{hi}$
that Strichartz won us to bound the terms appearing on the right-hand side of Proposition \ref{etprop}. For the rest
of this section all spacetime norms are going to be on $I_0\times\R^4$.

Let's consider \textbf{\eqref{et1}}. The term corresponding to $j=2$ can be bounded  via Bernstein and \eqref{u_hi estimate} by
$$
\eta_1 \|u_{hi}\|_{L^3_{t,x}}^3\|u_{lo}\|_{L^{\infty}_{t,x}} \lesssim \eta_1^3 C_0.
$$
We estimate the term corresponding to $j=1$ using Bernstein, \eqref{lowfreqcontrol}, and \eqref{u_hi estimate}:
$$
\eta_1 \|u_{hi}\|_{L^3_{t,x}}^2 \|u_{lo}\|_{L^{\infty}_tL^4_x} \|u_{lo}\|_{L^3_tL^{12}_x}
   \lesssim \eta_1^2 \eta_1^{\frac{7}{6}} C_{0}^{\frac{7}{6}}.
$$
The term coming from $j=0$ we again estimate using Bernstein and \eqref{lowfreqcontrol}:
\begin{align*}
\eta_1 \|u_{hi}P_{hi}\O(u_{lo}^3)\|_{L^1_{t,x}}
     &\lesssim \eta_1 \|u_{hi}\|_{L_t^{\infty}L_x^{2}} \|P_{hi}\O(u_{lo}^3)\|_{L_t^1L_x^2}
       \lesssim \eta_1^2 \|\nabla P_{hi}\O(u_{lo}^3)\|_{L_t^1L_x^2} \\
     &\lesssim \eta_1^2 \|\nabla u_{lo}\|_{L_t^2L_x^4} \|u_{lo}\|_{L^4_tL^8_x}^2
       \lesssim \eta_1^2 \|u_{lo}\|_{\dot{S}^1}^3
       \lesssim \eta_1^2 \eta_1^{\frac{3}{2}} C_0^{\frac{3}{2}}.
\end{align*}
The final contribution of the error term \eqref{et1} is therefore at most
$$\eta_1^2 (\eta_1 C_0 + \eta_1^{\frac{7}{6}} C_{0}^{\frac{7}{6}} + \eta_1^{\frac{3}{2}} C_0^{\frac{3}{2}}) \lesssim \eta_1.$$

We now consider the error term \textbf{\eqref{et2}}, which we can bound as
\begin{align*}
\eta_1 \int_{I_{0}} \int_{\R^4} |P_{lo}\O(u_{hi}^{3})(t,y)||u_{hi}(t,y)|dydt
    &\lesssim \eta_1 \|u_{hi}\|_{L_t^{\infty}L_x^2} \|P_{lo}\O(u_{hi}^{3})\|_{L_t^1L_x^2} \\
    &\lesssim \eta_1^2 \|P_{lo}\O(u_{hi}^{3})\|_{L^1_{t,x}}
      \lesssim \eta_1^2 \|u_{hi}\|_{L^3_{t,x}}^3\\
    &\lesssim \eta_1^3 C_0.
\end{align*}
Hence, \eqref{et2} is bounded by $\eta_1^3 C_0 \lesssim \eta_1$.

Consider \textbf{\eqref{et3}}. The contribution coming from $j=1$ can be estimated via \eqref{lowfreqcontrol} as
$$
\eta_1^2 \|\nabla u_{lo}\|_{L^2_tL^4_x} \|u_{lo}\|_{L^3_tL^{12}_x} \|u_{lo}\|_{L^6_{t,x}} \|u_{hi}\|_{L^{\infty}_tL^2_x}
   \lesssim \eta_1^3 \|u_{lo}\|_{\dot{S}^1}^3
   \lesssim \eta_1^3 \eta_1^{\frac{3}{2}} C_0^{\frac{3}{2}}.
$$
Applying Bernstein, \eqref{lowfreqcontrol}, and \eqref{u_hi estimate}, we estimate the contribution from $j=2$ by
$$
\eta_1 ^2 \|\nabla u_{lo}\|_{L^3_{t,x}} \|u_{hi}\|_{L^3_{t,x}}^2 \|u_{lo}\|_{L^{\infty}_{t,x}}
   \lesssim \eta_1^3 \eta_1^{\frac{7}{6}} C_0^{\frac{7}{6}}.
$$
Using Bernstein and \eqref{u_hi estimate} we estimate the term corresponding to $j=3$ by
$$
\eta_1^2 \|u_{hi}\|_{L^3_{t,x}}^3\|\nabla u_{lo}\|_{L^{\infty}_{t,x}}
    \lesssim \eta_1^3 C_0 \| u_{lo}\|_{L^{\infty}_tL^4_x}
    \lesssim \eta_1^4 C_0.
$$
Thus the error term \eqref{et3} contributes at most $\eta_1^3(\eta_1 C_0 + \eta_1^{\frac{7}{6}} C_{0}^{\frac{7}{6}}
 + \eta_1^{\frac{3}{2}} C_0^{\frac{3}{2}}) \lesssim \eta_1.$

We now turn towards  the error term \textbf{\eqref{et4}}. Decomposing $u=u_{hi}+u_{lo}$, we write
$\nabla P_{lo}(|u|^{2}u)=\sum_{j=0}^{3}\nabla P_{lo}\O(u_{hi}^{j}u_{lo}^{3-j})$.
Using Bernstein, Proposition \ref{Slf}, and Corollary \ref{cor uhi}, we estimate
\begin{align*}
\eta_1^2 \|u_{hi}\nabla P_{lo}\O(u_{lo}^3)\|_{L^1_{t,x}}
&\lesssim \eta_1^2 \|u_{hi}\|_{L_t^{\infty}L^2_x} \|\nabla u_{lo}\|_{L^2_tL^4_x} \|u_{lo}\|_{L^3_tL^{12}_x} \|u_{lo}\|_{L^6_{t,x}}
  \lesssim \eta_1^3 \|u_{lo}\|_{\dot{S}^1}^3\\
&\lesssim \eta_1^3 \eta_1^{\frac{3}{2}} C_0^{\frac{3}{2}}\lesssim \eta_1
\end{align*}
\begin{align*}
\eta_1^2 \|u_{hi}\nabla P_{lo}\O(u_{hi}u_{lo}^2)\|_{L^1_{t,x}}
&\lesssim \eta_1^2\|u_{hi}\|_{L_{t,x}^3}^2\|u_{lo}\|^2_{L_{t,x}^6}
  \lesssim \eta_1^2\|u_{hi}\|_{L_{t,x}^3}^2\|u_{lo}\|^2_{\dot{S}^1}\\
&\lesssim \eta_1^2 C_0^{2/3} \eta_1^{2/3} C_0 \eta_1\lesssim \eta_1\\
\eta_1^2 \|u_{hi}\nabla P_{lo}\O(u_{hi}^2u_{lo})\|_{L^1_{t,x}}
&\lesssim \eta_1^2 \|u_{hi}\|_{L^3_{t,x}}^3\|u_{lo}\|_{L_{t,x}^\infty}
  \lesssim \eta_1^2 \|u_{hi}\|_{L^3_{t,x}}^3\|u_{lo}\|_{L_t^\infty L_x^4}\\
&\lesssim \eta_1^2 C_0\eta_1 \eta_1\lesssim \eta_1
\end{align*}
\begin{align*}
\eta_1^2 \|u_{hi}\nabla P_{lo}\O(u_{hi}^3)\|_{L^1_{t,x}}
&\lesssim \eta_1^2 \|u_{hi}\|_{L_t^{\infty}L^2_x} \|\nabla P_{lo}\O(u_{hi}^3)\|_{L^1_tL^2_x}
  \lesssim \eta_1^3 \|u_{hi}^3\|_{L^1_{t,x}}\\
&\lesssim \eta_1^3 \|u_{hi}\|_{L^3_{t,x}}^3
\lesssim \eta_1^4 C_0\lesssim \eta_1.
\end{align*}
Hence, the error term \eqref{et4} is bounded by $\eta_1$.

We consider now \textbf{\eqref{et5}}. For the case $j=3$, we split the region of integration into $|x-y|\leq 1$ and $|x-y|>1$.
On the first region we use Cauchy-Schwartz and Young's inequality, remembering that now we are
convolving with an $L^2_x$ function, obtaining the bound
$$
\||u_{hi}|^2\|_{L_t^{\infty}L_x^2} \|u_{hi}^3 u_{lo}\|_{L^1_{t,x}}
   \lesssim \|u_{hi}\|_{L_t^{\infty}L_x^4}^{\frac{1}{2}} \|u_{hi}\|_{L^3_{t,x}}^3 \|u_{lo}\|_{L^{\infty}_{t,x}}
   \lesssim \eta_1^2 C_0.
$$
On the second region of integration, we bound $\frac{1}{|x-y|}<1$ and estimate
$$
\|u_{hi}\|_{L_t^{\infty}L_x^2}^2 \|u_{hi}^3 u_{lo}\|_{L^1_{t,x}}
   \lesssim \eta_1^2 \|u_{hi}\|_{L^3_{t,x}}^3 \|u_{lo}\|_{L^{\infty}_{t,x}}
   \lesssim \eta_1^4 C_0.
$$

The case $j=2$ is treated analogously using
$$
\|u_{hi}^2u_{lo}^2\|_{L^1_{t,x}}
   \lesssim \|u_{hi}\|_{L^{\infty}_tL^2_x} \|u_{hi}\|_{L^3_{t,x}} \|u_{lo}\|_{L^3_tL^{12}_x}^2
   \lesssim \eta_1 \eta_1^{\frac{1}{3}} C_0^{\frac{1}{3}} \|u_{lo}\|_{\dot{S}^1}^2
   \lesssim \eta_1^{\frac{7}{3}} C_0^{\frac{4}{3}}
$$
to get the bound $\eta_1^{\frac{7}{3}}C_0^{\frac{4}{3}}$ on the first region
and $\eta_1^2 \eta_1^{\frac{7}{3}}C_0^{\frac{4}{3}}$ on the second region.

We turn now towards the term corresponding to $j=1$. We write $u_{hi}=\nabla (\nabla^{-1}u_{hi})$
and integrate by parts to bound this term by
\begin{align}
&\Bigl|\int_{I_{0}}\int_{\R^4} \int_{\R^4} \frac{|u_{hi}(t,y)|^{2}(\nabla^{-1}u_{hi})(t,x)u_{lo}^2(t,x)\nabla u_{lo}(t,x)}{|x-y|}dxdydt\Bigr| \label{a1}\\
&+\Bigl|\int_{I_{0}}\int_{\R^4} \int_{\R^4} |u_{hi}(t,y)|^{2}u_{lo}^3(t,x)(\nabla^{-1}u_{hi})(t,x)\frac{x-y}{|x-y|^3}dxdydt\Bigr|.\label{a2}
\end{align}

Let's consider \eqref{a1}. We again split the domain of integration into $|x-y|\leq 1$ and $|x-y|>1$.
On the first region we proceed as in the previous two cases, i.e. use Cauchy-Schwartz, Young's inequality,
and that the function we convolve with is in $L^1_x$ to bound
\begin{align*}
\||u_{hi}|^2\|_{L_t^{\infty}L_x^2}&\|(\nabla^{-1}u_{hi})u_{lo}^2 \nabla u_{lo}\|_{L^1_tL^2_x}\\
    &\lesssim \|u_{hi}\|_{L_t^{\infty}L_x^4}^{\frac{1}{2}}\|\nabla^{-1}u_{hi}\|_{L^2_tL^4_x}
    \|\nabla u_{lo}\|_{L^2_tL^4_x}\|u_{lo}\|_{L^{\infty}_{t,x}}^2 \\
&\lesssim \eta_1^2 \|u_{lo}\|_{\dot{S}^1}\|\nabla^{-\frac{1}{2}}u_{hi}\|_{L^2_tL^4_x} \\
&\lesssim\eta_1^3 C_0.
\end{align*}
On the second region we take the absolute values inside the integral and estimate
\begin{align*}
\int_{I_{0}} \int_{\R^4} \int_{\R^4} &\frac{|u_{hi}(t,y)|^{2}|(\nabla^{-1}u_{hi})(t,x)||u_{lo}^2(t,x)||\nabla u_{lo}(t,x)|}{\langle x-y\rangle}dxdydt \\
                       &\lesssim \|u_{hi}\|_{L_t^{\infty}L_x^2}^2 \sup_{y\in\R^4} \int_{I_0} \int_{\R^4} \frac{|(\nabla^{-1}u_{hi})(t,x)||u_{lo}^2(t,x)||\nabla u_{lo}(t,x)|}{\langle x-y\rangle}dxdt \\
                       &\lesssim \eta_1^2 \|{\langle x\rangle}^{-1}\|_{L_x^{4+}} \|\nabla^{-1}u_{hi}u_{lo}^2\nabla u_{lo}\|_{L^1_tL^{4/3-}_x} \\
                       &\lesssim \eta_1^2 \|\nabla^{-1}u_{hi}\|_{L^2_tL^4_x} \|\nabla u_{lo}\|_{L^2_tL^4_x} \|u_{lo}\|_{L_t^{\infty}L_x^{8-}}^2 \\
                       &\lesssim \eta_1^4 \|u_{lo}\|_{\dot{S}^1} \|\nabla^{-\frac{1}{2}}u_{hi}\|_{L^2_tL^4_x} \\
                       &\lesssim \eta_1^5 C_0.
\end{align*}

Taking absolute values inside the integrals, we write \eqref{a2} as
\begin{align*}
|\langle\nabla^{-2}|u_{hi}|^2,|(\nabla^{-1}u_{hi})u_{lo}^3|\rangle|
   &= |\langle\nabla^{-\frac{1}{2}}|u_{hi}|^2,\nabla^{-\frac{3}{2}}|(\nabla^{-1}u_{hi})u_{lo}^3|\rangle| \\
   &\lesssim \|\nabla^{-\frac{1}{2}}|u_{hi}|^2\|_{L^2_{t,x}} \|\nabla^{-\frac{3}{2}}|(\nabla^{-1}u_{hi})u_{lo}^3|\|_{L^2_{t,x}} \\
   &\lesssim \eta_1^{\frac{1}{2}} C_0^{\frac{1}{2}} \Bigl\|\int_{\R^4}\frac{|\nabla^{-1}u_{hi}(t,x)||u_{lo}(t,x)|^3}{|x-y|^{\frac{5}{2}}}dx\Bigr\|_{L^2_{t,x}}.
\end{align*}
To estimate this integral we again split the domain of integration into $|x-y|\leq 1$ and $|x-y|>1$.
On the first region we note that the function we convolve with is in $L^1_x$ and we estimate
$$
\|(\nabla^{-1}u_{hi})u_{lo}^3\|_{L^2_{t,x}}
   \lesssim \|\nabla^{-1}u_{hi}\|_{L^2_tL^4_x} \|u_{lo}\|_{L^{\infty}_tL_x^8}^2 \|u_{lo}\|_{L^{\infty}_{t,x}}
   \lesssim \eta_1^3 \eta_1^{\frac{1}{2}} C_0^{\frac{1}{2}} \lesssim \eta_1^{\frac{7}{2}} C_0^{\frac{1}{2}}.
$$
On the second region of integration we estimate
\begin{align*}
\|\{(\nabla^{-1}u_{hi})u_{lo}^3\} * \langle x\rangle^{-\frac{5}{2}}\|_{L^2_{t,x}}
  &\lesssim \|(\nabla^{-1}u_{hi})u_{lo}^3\|_{L^2_tL^1_x}\|\langle x\rangle^{-\frac{5}{2}}\|_{L^{\infty}_tL^2_x} \\
  &\lesssim\|\nabla^{-1}u_{hi}\|_{L^2_tL^4_x}\|u_{lo}\|_{L^{\infty}_tL_x^4}^3 \\
  &\lesssim \eta_1^3 \eta_1^{\frac{1}{2}}C_0^{\frac{1}{2}}
  \lesssim \eta_1^{\frac{7}{2}}C_0^{\frac{1}{2}}.
\end{align*}
We see therefore that the error term \eqref{et5} is also bounded by $\eta_1$.

Upon rescaling, this concludes the proof of Proposition \ref{prop flim}.

%%%%%%%%%%%%%%%%%%%%%%%%%%%%%%%%%%%%%%%%%%%%%%%%%%%%%%%%%%%%%%%%%%%%%%%%%%%%%%%%%%%%%%%%%%%
%
%
%                                   Section
%
%
%%%%%%%%%%%%%%%%%%%%%%%%%%%%%%%%%%%%%%%%%%%%%%%%%%%%%%%%%%%%%%%%%%%%%%%%%%%%%%%%%%%%%%%%%%%

\section{Preventing Energy Evacuation}
The purpose of this section is to prove

\begin{proposition}[Energy cannot evacuate to high frequencies]\label{energyevac} Suppose that $u$ is a minimal energy blowup solution
of \eqref{schrodinger equation}. Then for all $t\in I_{0}$,
\begin{align}
N(t)\lesssim C(\eta_3)N_{min}. \label{freqbound}
\end{align}
\end{proposition}

\begin{remark}\label{rem Nmin}
$N_{min}>0$. Indeed, from the frequency localization result we know that for $\forall t \in I_0$
$$
\|P_{c(\eta_0)N(t)<\cdot<C(\eta_0)N(t)}u(t)\|_{\dot{H}^1_x}\sim 1.
$$
On the other hand, from Bernstein
$$
\|P_{c(\eta_0)N(t)<\cdot<C(\eta_0)N(t)}u(t)\|_{\dot{H}^1_x}
   \leq C(\eta_0) N(t) \|u\|_{L_t^{\infty}L_x^2}.
$$
Hence, $N(t) \geq c(\eta_0) \|u\|_{L_t^{\infty}L_x^2}^{-1}$ for all $t \in I_0$,
which implies $N_{min}=\inf_{t \in I_0} N(t) > 0$.
\end{remark}

\subsection{The setup}
We normalize $N_{min}=1$. As $N(t) \in 2^{\Z}$, there exists $t_{min} \in I_0$ such that $N(t_{min})=N_{min}=1$.

At the time $t=t_{min}$ we have a considerable amount of mass at medium frequencies:
\begin{align}
\|P_{c(\eta_0)<\cdot<C(\eta_0)}u(t_{min})\|_{L^2_x}
   \gtrsim c(\eta_0)\|P_{c(\eta_0)<\cdot<C(\eta_0)}u(t_{min})\|_{\dot{H}^1_x}
   \sim c(\eta_0). \label{massmedfreq}
\end{align}
However, by Bernstein there is not much mass at frequencies higher than $C(\eta_0)$:
$$
\|P_{>C(\eta_0)}u(t_{min})\|_{L^2_x} \lesssim c(\eta_0).
$$

Let's assume for a contradiction that there exists $t_{evac} \in I_0$ such that $N(t_{evac})\gg C(\eta_3)$.
By time reversal symmetry we may assume $t_{evac} > t_{min}$. As
$$\|P_{<c(\eta)N(t)}u(t)\|_{\dot{H}^1_x} \leq \eta$$ for every $\eta_3 \leq \eta \leq \eta_0$ and all $t \in I_0$,
we see that choosing $C(\eta_3)$ sufficiently large we get very small energy at low and medium frequencies at the time $t=t_{evac}$:
\begin{align}
\|P_{< \eta_3^{-1}}u(t_{evac})\|_{\dot{H}^1} \leq \eta_3. \label{goodcontrollowfreq}
\end{align}

We define $u_{lo}=P_{<\eta_{2}^{100}}u$ and $u_{hi}=P_{\geq \eta_{2}^{100}}u$. \eqref{massmedfreq} implies that
\begin{align}
\|u_{hi}(t_{min})\|_{L^2_x} \geq \eta_1. \label{highmass}
\end{align}
Suppose we could show that a big portion of the mass sticks around until $t=t_{evac}$, i.e.
\begin{align}
\|u_{hi}(t_{evac})\|_{L^2_x} \geq \frac{1}{2}\eta_1. \label{highmassasmp}
\end{align}
Then, since by Bernstein
$$
\|P_{>C(\eta_1)}u_{hi}(t_{evac})\|_{L^2_x}\leq c(\eta_1),
$$
the triangle inequality would imply
$$
\|P_{\leq C(\eta_1)}u_{hi}(t_{evac})\|_{L^2_x} \geq \frac{1}{4}\eta_1.
$$
Another application of Bernstein yields
$$
\|P_{\leq C(\eta_1)}u(t_{evac})\|_{\dot{H}^1_x} \gtrsim c(\eta_1, \eta_2),
$$
which would contradict \eqref{goodcontrollowfreq} if $\eta_3$ were chosen sufficiently small.

It therefore remains to show \eqref{highmassasmp}. In order to prove it we assume that there exists a time $t_{*}$
such that $t_{min} \leq t_{*} \leq t_{evac}$ and
\begin{align}
\inf_{t_{min} \leq t \leq t_{*}}\|u_{hi}(t)\|_{L^2_x} \geq \frac{1}{2} \eta_1. \label{contasmp}
\end{align}
We will show that this can be bootstrapped to
\begin{align}
\inf_{t_{min} \leq t \leq t_{*}}\|u_{hi}(t)\|_{L^2_x} \geq \frac{3}{4} \eta_1. \label{bootstrapmass}
\end{align}
Hence $\{t_{*} \in [t_{min}, t_{evac}] : \eqref{contasmp} \text{ holds} \}$ is both open and closed in $[t_{min}, t_{evac}]$,
and \eqref{highmassasmp} holds.

In order to show that \eqref{contasmp} implies \eqref{bootstrapmass} we treat the $L^2_x$-norm of $u_{hi}$
as an almost conserved quantity. Define
$$
L(t)=\int_{\R^4}|u_{hi}(t,x)|^2dx.
$$
By \eqref{highmass} we have $L(t_{min}) \geq \eta_1^2$. Hence,
by the Fundamental Theorem of Calculus it suffices to show that
$$
\int_{t_{min}}^{t_{*}}|\partial_{t}L(t)|dt \leq \frac{1}{100}\eta_1^2.
$$

We have
\begin{align*}
\partial_{t}L(t)
  &=2 \int_{\R^4} \{P_{hi}(|u|^2u), u_{hi}\}_{m}dx \\
  &=2 \int_{\R^4} \{P_{hi}(|u|^2u)-|u_{hi}|^2u_{hi}, u_{hi}\}_{m}dx.
\end{align*}
Thus, we need to show
\begin{align}
\int_{t_{min}}^{t_{*}} \Bigl| \int_{\R^4} \{P_{hi}(|u|^2u)-|u_{hi}|^2u_{hi}, u_{hi}\}_{m}dx\Bigr| dt \leq \frac{1}{100}\eta_1^2. \label{last}
\end{align}

In order to prove \eqref{last}, we need to control the various interactions between low, medium, and high frequencies.
In the next section we will develop the tools that will make this goal possible.

\subsection{Spacetime estimates on low, medium, and high frequencies}

Remember that the frequency-localized Morawetz inequality implies that for $N<c(\eta_1)$,
\begin{align}
\int_{t_{min}}^{t_{evac}} \int_{\R^4} |P_{\geq N}u(t,x)|^3 dxdt \lesssim \eta_1 N^{-3}. \label{conseqflim}
\end{align}

This estimate is useful for medium and high frequencies; however it is extremely bad for low frequencies since
$N^{-3}$ diverges as $N \rightarrow 0$. We therefore need to develop better estimates in this case.
As $u_{\leq \eta_2}$ has extremely small energy at $t=t_{evac}$ (see \eqref{goodcontrollowfreq}),
we expect to have small energy at all times in $[t_{min}, t_{evac}]$. Of course, there is energy leaking from
the high frequencies to the low frequencies, but \eqref{conseqflim} limits this leakage. Indeed, we have

\begin{lemma}
Under the assumptions above, we have
\begin{align}
\|P_{\leq N}u\|_{\dot{S}^1([t_{min}, t_{evac}]\times \R^4)}
    \lesssim \eta_3 + \eta_2^{-\frac{3}{2}} N^{\frac{3}{2}} \label{S^1lf}
\end{align}
for all $N \leq \eta_2$.
\end{lemma}

\begin{remark}
One should think of the $\eta_3$ factor on the right-hand side of \eqref{S^1lf} as the energy coming from
the low modes of $u(t_{evac})$, and the $\eta_2^{-\frac{3}{2}} N^{\frac{3}{2}}$ term as the energy coming from the
high frequencies of $u(t)$ for $t_{min} \leq t \leq t_{evac}$.
\end{remark}

\begin{proof}
Consider the set
$$
\Omega=\{t \in [t_{min}, t_{evac}) : \|P_{\leq N}u\|_{\dot{S}^1([t, t_{evac}]\times \R^4)}
       \leq C_0 \eta_3 + \eta_0 \eta_2^{-\frac{3}{2}} N^{\frac{3}{2}}, \forall N \leq \eta_2\},
$$
where $C_0$ is a large constant to be chosen later and not depending on any of the $\eta$'s.

Our goal is to show that $t_{min} \in \Omega$. First we will show that $t \in \Omega$ for $t$ close to $t_{evac}$.
Indeed, from Strichartz and Sobolev  we get
\begin{align*}
\|P_{\leq N}u\|_{\dot{S}^1([t, t_{evac}]\times \R^4)}
    \lesssim & \|\nabla P_{\leq N}u\|_{L_t^{\infty}L_x^2([t, t_{evac}]\times \R^4)}
               + \|\nabla u\|_{L^2_tL^4_x([t, t_{evac}]\times \R^4)} \\
    \lesssim & \|\nabla P_{\leq N}u(t_{evac})\|_{L_x^2} + C|t_{evac}-t| \|\nabla \partial_{t}u\|_{L_t^{\infty}L_x^2(I_0\times \R^4)} \\
             & + |t_{evac}-t|^{\frac{1}{2}} \|\nabla u\|_{L^{\infty}_tL^4_x(I_0\times \R^4)}.
\end{align*}
As $u$ is Schwartz the last two norms are finite, so \eqref{goodcontrollowfreq} implies
$$
\|P_{\leq N}u\|_{\dot{S}^1([t, t_{evac}]\times \R^4)}
   \lesssim \eta_3 + C(I_0,u) |t_{evac}-t| + C(I_0,u)|t_{evac}-t|^{\frac{1}{2}}.
$$
Thus $t \in \Omega$ provided $|t_{evac}-t|$ is sufficiently small and we choose $C_0$ sufficiently large.

Now suppose that $t \in \Omega$. We will show that
\begin{align}
\|P_{\leq N}u\|_{\dot{S}^1([t, t_{evac}]\times \R^4)}
       \leq \frac{1}{2} C_0 \eta_3 + \frac{1}{2} \eta_0 \eta_2^{-\frac{3}{2}} N^{\frac{3}{2}} \label{bootstraplast}
\end{align}
holds for $\forall N \leq \eta_2$. Thus, $\Omega$ is both open and closed in $[t_{min},t_{evac}]$ and we have $t_{min} \in \Omega$ as desired.

Fixing $N \leq \eta_2$, Lemma \ref{lemma linear strichartz} implies
$$
\|P_{\leq N}u\|_{\dot{S}^1([t, t_{evac}]\times\R^4)}
   \leq \|P_{\leq N}u(t_{evac})\|_{\dot{H}^1_x} + \sum_{m=1}^{M} \|\nabla F_{m}\|_{L_t^{q'_m}L_x^{r'_m}([t, t_{evac}]\times\R^4)}
$$
for some decomposition $P_{\leq N}(|u|^2u)=\sum_{m=1}^M F_m$ and some dual Schr\"odinger admissible pair $(q'_m,r'_m)$.

From \eqref{goodcontrollowfreq} we have
$$
\|P_{\leq N}u(t_{evac})\|_{\dot{H}^1_x} \lesssim \eta_3,
$$
which is acceptable if we choose $C_0$ sufficiently large.

We decompose $P_{\leq N}(|u|^2u)=\sum_{j=0}^3 P_{\leq N}\O(u_{hi'}^ju_{lo'}^{3-j}) =\sum_{j=0}^3 F_j$ with
$u_{lo'}=u_{\leq \eta_2}$ and $u_{hi'}=u_{>\eta_2}$.

Consider the $j=3$ term. Using Bernstein and \eqref{goodlemma} we estimate
\begin{align*}
\|\nabla P_{\leq N}\O(u_{hi'}^3)\|_{L^2_tL_x^{4/3}([t, t_{evac}]\times\R^4)}
   &\lesssim N^{\frac{3}{2}} \|\nabla^{-\frac{1}{2}}P_{\leq N}\O(u_{hi'}^3)\|_{L^2_tL_x^{4/3}([t, t_{evac}]\times\R^4)} \\
   &\lesssim N^{\frac{3}{2}} (\eta_1 \eta_2^{-3})^{\frac{1}{2}}
    =\eta_1^{\frac{1}{2}} \eta_2^{-\frac{3}{2}} N^{\frac{3}{2}}.
\end{align*}

Now consider the term corresponding to $j=2$. By Bernstein and H\"older,
\begin{align*}
\|\nabla P_{\leq N}\O(&u_{hi'}^2 u_{lo'})\|_{L^2_tL_x^{4/3}([t, t_{evac}]\times\R^4)}\\
   &\lesssim N^2 \|u_{hi'}^2 u_{lo'}\|_{L^2_tL_x^1([t, t_{evac}]\times\R^4)} \\
   &\lesssim N^2 \|u_{hi'}\|_{L^3_{t,x}([t, t_{evac}]\times\R^4)} \|u_{hi'}\|_{L_t^{\infty}L_x^2([t, t_{evac}]\times\R^4)} \|u_{lo'}\|_{L^6_{t,x}([t, t_{evac}]\times\R^4)} \\
   &\lesssim N^2 (\eta_1 \eta_2^{-3})^{\frac{1}{3}} \eta_2^{-1} \|u_{lo'}\|_{\dot{S}^1([t, t_{evac}]\times\R^4)}.
\end{align*}
But since $t \in \Omega$, $\|u_{lo'}\|_{\dot{S}^1([t, t_{evac}]\times\R^4)}\leq \eta_0$, and
$$
\|\nabla P_{\leq N}\O(u_{hi'}^2 u_{lo'})\|_{L^2_tL_x^{4/3}([t, t_{evac}]\times\R^4)}
    \lesssim \eta_0 \eta_2^{-2} N^2
    \lesssim \eta_0 \eta_2^{-\frac{3}{2}} N^{\frac{3}{2}}.
$$

We turn now towards the $j=1$ contribution.  Consider first the case when $N<c\eta_2$. The expression $\|\nabla P_{\leq N}\O(u_{hi'} u_{lo'}^2)\|_{L^{q'_1}_tL_x^{r'_1}([t, t_{evac}]\times\R^4)}$
vanishes unless one of the $u_{lo'}$ factors has frequency $>c\eta_2$. Writing $u_{lo'}=P_{\leq c\eta_2}u_{lo'}+P_{> c\eta_2}u_{lo'}$
and recalling that the projections are bounded on the space considered, we see that we only need to control
$$
\|\nabla P_{\leq N}\O(u_{hi'} (P_{>c\eta_2}u_{lo'}) u_{lo'})\|_{L^2_tL_x^{4/3}([t, t_{evac}]\times\R^4)}.
$$
But $P_{>c\eta_2}u_{lo'}$ obeys the same $L^3_{t,x}$ estimates as $u_{hi'}$, so controlling the above expression
amounts to the manipulations of case $j=2$.

Now consider the case when $N \geq c\eta_2$. Taking $(q'_1, r'_1)=(1, 2)$ and using Bernstein we have
\begin{align*}
\|\nabla P_{\leq N}\O(u_{hi'} u_{lo'}^2)&\|_{L^1_tL_x^2([t, t_{evac}]\times\R^4)}\\
   &\lesssim \eta_2 \|u_{hi'} u_{lo'}^2\|_{L^1_tL_x^2([t, t_{evac}]\times\R^4)} \\
   &\lesssim \eta_2 \|u_{hi'}\|_{L^3_{t,x}([t, t_{evac}]\times\R^4)} \|u_{lo'}\|_{L^3_tL^{12}_x([t, t_{evac}]\times\R^4)}^2 \\
   &\lesssim \eta_2 (\eta_1 \eta_2^{-3})^{\frac{1}{3}} \|u_{lo'}\|_{\dot{S}^1([t, t_{evac}]\times\R^4)}^2 \\
   &\lesssim \eta_0^2 \eta_1^{\frac{1}{3}}
    \lesssim \eta_0 \eta_2^{-\frac{3}{2}} N^{\frac{3}{2}}.
\end{align*}

We are left with the $j=0$ term. Write $u_{lo'}=u_{<\eta_3}+u_{\eta_3\leq \cdot \leq \eta_2}$. Any term containing at least
one $u_{<\eta_3}$ can be controlled using the trilinear Strichartz estimates
\begin{align*}
\|\nabla P_{\leq N}\O(u_{lo'}^2u_{<\eta_3})&\|_{L_t^1L_x^2([t, t_{evac}]\times\R^4)}\\
   &\lesssim \|u_{lo'}\|_{\dot{S}^1([t, t_{evac}]\times\R^4)}^2 \|u_{<\eta_3}\|_{\dot{S}^1([t, t_{evac}]\times\R^4)} \\
   &\lesssim (C_0 \eta_3 + \eta_0)^2 (C_0 \eta_3 + \eta_0 \eta_2^{-\frac{3}{2}} \eta_3^{\frac{3}{2}}) \\
   &\lesssim C_0 \eta_0^2 \eta_3,
\end{align*}
which is acceptable if we choose $\eta_3$ sufficiently small.

We can thus discard all the terms involving $u_{<\eta_3}$ and focus on the term
$$\|\nabla P_{\leq N} \O(u_{\eta_3 \leq \cdot \leq \eta_2}^3)\|_{L_t^1L_x^2([t, t_{evac}]\times\R^4)}.$$
Using Bernstein we estimate this as
\begin{align*}
\|\nabla P_{\leq N} \O(u_{\eta_3 \leq \cdot \leq \eta_2}^3)\|_{L_t^1L_x^2([t, t_{evac}]\times\R^4)}
   &\lesssim N \| P_{\leq N} \O(u_{\eta_3 \leq \cdot \leq \eta_2}^3)\|_{L_t^1L_x^2([t, t_{evac}]\times\R^4)} \\
   &\lesssim N^{\frac{3}{2}} \|u_{\eta_3 \leq \cdot \leq \eta_2}^3\|_{L_t^1L_x^{8/5}([t, t_{evac}]\times\R^4)} \\
   &\lesssim N^{\frac{3}{2}} \|u_{\eta_3 \leq \cdot \leq \eta_2}\|_{L_t^3L_x^{24/5}([t, t_{evac}]\times\R^4)}^3.
\end{align*}
But by Bernstein and the hypothesis $t\in\Omega$, we have
\begin{align*}
\|u_{\eta_3 \leq \cdot \leq \eta_2}\|_{L_t^3L_x^{24/5}([t, t_{evac}]\times\R^4)}
   &\lesssim \sum_{\eta_3 \leq M \leq \eta_2} \|P_M u\|_{L_t^3L_x^{24/5}([t, t_{evac}]\times\R^4)} \\
   &\lesssim \sum_{\eta_3 \leq M \leq \eta_2} M^{-1} \|\nabla P_M u\|_{L_t^3L_x^{24/5}([t, t_{evac}]\times\R^4)} \\
   &\lesssim \sum_{\eta_3 \leq M \leq \eta_2} M^{-\frac{1}{2}} \|\nabla P_M u\|_{L^3_{t,x}([t, t_{evac}]\times\R^4)} \\
   &\lesssim \sum_{\eta_3 \leq M \leq \eta_2} M^{-\frac{1}{2}} \| P_M u\|_{\dot{S}^1([t, t_{evac}]\times\R^4)} \\
   &\lesssim \sum_{\eta_3 \leq M \leq \eta_2} M^{-\frac{1}{2}} (C_0 \eta_3 + \eta_0 \eta_2^{-\frac{3}{2}} M^{\frac{3}{2}}) \\
   &\lesssim \eta_0 \eta_2^{-\frac{1}{2}}.
\end{align*}
Hence,
$$
\|\nabla P_{\leq N} \O(u_{\eta_3 \leq \cdot \leq \eta_2}^3)\|_{L_t^1L_x^2([t, t_{evac}]\times\R^4)}
    \lesssim \eta_0^3 \eta_2^{-\frac{3}{2}} N^{\frac{3}{2}}.
$$
This proves \eqref{bootstraplast} and closes the bootstrap.

\end{proof}

\subsection{Controlling the localized mass increment.} We now have good enough control over low, medium, and high
frequencies to prove \eqref{last}.  Writing
\begin{align*}
P_{hi}(|u|^2u)-|u_{hi}|^2u_{hi}
  &=P_{hi}(|u|^2u-|u_{hi}|^2u_{hi}-|u_{lo}|^2u_{lo}) + P_{hi}(|u_{lo}|^2u_{lo})\\
  &\quad - P_{lo}(|u_{hi}|^2u_{hi}),
\end{align*}
we only have to consider the terms
\begin{align}
\int_{t_{min}}^{t_{*}}\Bigl|\int_{\R^4} \overline{u_{hi}}P_{hi}(|u|^2u-|u_{hi}|^2u_{hi}-|u_{lo}|^2u_{lo})dx\Bigr|dt \label{1} \\
\int_{t_{min}}^{t_{*}}\Bigl|\int_{\R^4} \overline{u_{hi}}P_{hi}(|u_{lo}|^2u_{lo})dx\Bigr|dt \label{2} \\
\int_{t_{min}}^{t_{*}}\Bigl|\int_{\R^4} \overline{u_{hi}}P_{lo}(|u_{hi}|^2u_{hi})dx\Bigr|dt \label{3}.
\end{align}

Take \eqref{1}. We use the inequality
$$
||u|^2u-|u_{hi}|^2u_{hi}-|u_{lo}|^2u_{lo}|\lesssim |u_{hi}|^2|u_{lo}| + |u_{hi}||u_{lo}|^2
$$
to estimate
\begin{align}
\eqref{1} &\lesssim \int_{t_{min}}^{t_{*}}\int_{\R^4} |P_{hi}u_{hi}|(|u_{hi}||u_{lo}|^2 + |u_{hi}|^2|u_{lo}|)dxdt \\
          &= \int_{t_{min}}^{t_{*}}\int_{\R^4} |P_{hi}u_{hi}||u_{hi}||u_{lo}|^2 dxdt \label{ts1} \\
          &\quad+ \int_{t_{min}}^{t_{*}}\int_{\R^4} |P_{hi}u_{hi}||u_{hi}|^2|u_{lo}| dxdt \label{ts2}.
\end{align}
We estimate
\begin{align*}
\eqref{ts1}
   &\lesssim \|P_{hi}u_{hi}\|_{L^3_{t,x}([t_{min}, t_{*}]\times\R^4)} \|u_{hi}\|_{L^3_{t,x}([t_{min}, t_{*}]\times\R^4)} \|u_{lo}\|_{L^6_{t,x}([t_{min}, t_{*}]\times\R^4)}^2 \\
   &\lesssim \{\eta_1 (\eta_2^{100})^{-3}\}^{\frac{2}{3}} \{\eta_3 + \eta_2^{-\frac{3}{2}} (\eta_2^{100})^{\frac{3}{2}}\}^2 \\
   &\lesssim \eta_1^{\frac{2}{3}} \eta_2^{97} \ll \eta_1^2.
\end{align*}
We now turn towards the contribution of \eqref{ts2}. We decompose
$u_{hi}=u_{\eta_2^{100} \leq \cdot \leq \eta_2} + u_{>\eta_2}$ and estimate
$$
\|u_{>\eta_2}\|_{L^3_{t,x}([t_{min}, t_{*}]\times\R^4)} \lesssim (\eta_1 \eta_2^{-3})^{\frac{1}{3}}= \eta_1^{\frac{1}{3}} \eta_2^{-1}
$$
and
\begin{align*}
\|u_{\eta_2^{100} \leq \cdot \leq \eta_2}\|_{L^3_{t,x}([t_{min}, t_{*}]\times\R^4)}
   &\lesssim \sum_{\eta_2^{100} \leq N \leq \eta_2} \|u_N\|_{L^3_{t,x}([t_{min}, t_{*}]\times\R^4)} \\
   &\lesssim \sum_{\eta_2^{100} \leq N \leq \eta_2} N^{-1} \|\nabla u_N\|_{L^3_{t,x}([t_{min}, t_{*}]\times\R^4)} \\
   &\lesssim \sum_{\eta_2^{100} \leq N \leq \eta_2} N^{-1} \|u_N\|_{\dot{S}^1([t_{min}, t_{*}]\times\R^4)} \\
   &\lesssim \sum_{\eta_2^{100} \leq N \leq \eta_2} N^{-1} (\eta_3 + \eta_2^{-\frac{3}{2}} N^{\frac{3}{2}}) \\
   &\lesssim \eta_2^{-\frac{3}{2}} \eta_2^{\frac{1}{2}}
    =\eta_2^{-1}.
\end{align*}
By Bernstein and \eqref{potential energy bound},
$$\|u_{lo}\|_{L^{\infty}_{t,x}([t_{min}, t_{*}]\times\R^4)} \lesssim \eta_2^{100} \|u_{lo}\|_{L^{\infty}_tL_x^4([t_{min}, t_{*}]\times\R^4)} \lesssim \eta_2^{100}.$$
As $P_{hi}u_{hi}$ obeys the same estimates as $u_{hi}$, we bound the contribution of \eqref{ts2} by
\begin{align*}
&\int_{t_{min}}^{t_{*}}\int_{\R^4} |u_{hi}|^3|u_{lo}|dxdt\\
   &\quad \lesssim \sum_{j=0}^3 \|u_{\eta_2^{100} \leq \cdot \leq \eta_2}\|_{L^3_{t,x}([t_{min}, t_{*}]\times\R^4)}^j
      \|u_{>\eta_2}\|_{L^3_{t,x}([t_{min}, t_{*}]\times\R^4)}^{3-j} \|u_{lo}\|_{L^{\infty}_{t,x}([t_{min}, t_{*}]\times\R^4)}\\
   &\quad\lesssim \eta_2^{-3} \eta_2^{100} \ll \eta_1^2.
\end{align*}
Hence, $\eqref{1} \ll \eta_1^2.$

We consider next \eqref{2}. Because of the presence of $P_{hi}$, one of the terms $u_{lo}$ must have frequency larger than
$c\eta_3$ or the expression vanishes. Moving $P_{hi}$ over to $\overline{u_{hi}}$, we bound \eqref{2} by a sum of terms that look like
$$
\int_{t_{min}}^{t_{*}}\int_{\R^4} |P_{hi}u_{hi}||P_{\geq c\eta_2^{100}}u_{lo}||u_{lo}|^2dxdt.
$$
As $P_{\geq c\eta_2^{100}}u_{lo}=P_{lo}u_{\geq c\eta_2^{100}}$ behaves like $u_{hi}$ as far as norms are concerned,
this expression can be estimated by the procedure used to estimate \eqref{ts1}.

We now turn to \eqref{3}. Moving the projection $P_{lo}$ onto $\overline{u_{hi}}$ and writing
$P_{lo}u_{hi}=P_{hi}u_{lo}$, we bound the contribution of \eqref{3} by
$$
\int_{t_{min}}^{t_{*}}\int_{\R^4} |P_{hi}u_{lo}| |u_{hi}|^3 dxdt.
$$
As $P_{hi}u_{lo}$ behaves like $u_{lo}$ as far as norms are concerned, \eqref{3} is bounded by the same method as
\eqref{ts2}.

Therefore \eqref{last} holds, and this concludes the proof of Proposition \ref{energyevac}.

%%%%%%%%%%%%%%%%%%%%%%%%%%%%%%%%%%%%%%%%%%%%%%%%%%%%%%%%%%%%%%%%%%%%%%%%%%%%%%%%%%%%%%%%%%%
%
%
%                                   Section
%
%
%%%%%%%%%%%%%%%%%%%%%%%%%%%%%%%%%%%%%%%%%%%%%%%%%%%%%%%%%%%%%%%%%%%%%%%%%%%%%%%%%%%%%%%%%%%

\section{The contradiction argument}

We now have all the information we need about a minimal energy blowup solution to conclude the
contradiction argument. We know it is localized in frequency and
concentrates in space. The Morawetz inequality provides good control on the high-frequency part
of $u$ in $L^3_{t,x}$.  By the arguments in the
previous section we have excluded the last enemy by showing that the solution can't shift its energy from low modes to high
modes causing the $L^6_{t,x}$-norm to blow up while the $L^3_{t,x}$-norm remains bounded. Hence $N(t)$ must remain within a
bounded set $[N_{min},N_{max}]$, where $N_{max}\leq C(\eta_3)N_{min}$ and
$N_{min}\geq c(\eta_0)\|u\|_{L^{\infty}_tL_x^2}^{-1}$.
Combining all these (and again relying on the Morawetz inequality) we can derive the desired contradiction.
We begin with

\begin{lemma}\label{lemma integral bound on N} For any minimal
energy blowup solution of \eqref{schrodinger equation} we have
\begin{equation}\label{integral bound on N}
\int_{I_0} N(t)^{-1} dt \lesssim C(\eta_1,\eta_2) N_{min}^{-3}.
\end{equation}
In particular, because $N(t) \leq C(\eta_3) N_{min}$ for all $t \in I_0$, we have
\begin{equation}\label{length I bound}
|I_0| \lesssim C(\eta_1, \eta_2, \eta_3,N_{min}).
\end{equation}
\end{lemma}

\begin{proof}
By \eqref{u_hi estimate} we have
\begin{equation*}
\int_{I_0} \int_{\R^4} |P_{\geq N_*}u|^3 dx dt \lesssim \eta_1
N_*^{-3}
\end{equation*}
for all $N_* < c(\eta_1)N_{min}$.  Let $N_* = c(\eta_2) N_{min}$ and rewrite the above estimate as
\begin{equation}\label{int bound low}
\int_{I_0} \int_{\R^4} |P_{\geq N_*}u|^3 dx dt \lesssim C(\eta_1,\eta_2)N_{min}^{-3}.
\end{equation}

On the other hand, by Bernstein and \eqref{potential energy bound}
we have
\begin{equation}\label{int bound high}
\begin{split}
\int_{|x-x(t)| \leq C(\eta_1)/N(t)}|P_{<N_*}u(t)|^3 dx & \lesssim
C(\eta_1)N(t)^{-4} \|P_{<N_*} u(t)\|_{L^\infty_x}^3\\
& \lesssim C(\eta_1)N(t)^{-4}N(t)^3 c(\eta_2)\|P_{<N_*}u(t)\|_{L^4_x}^3\\
& \lesssim c(\eta_2) N(t)^{-1}.
\end{split}
\end{equation}
By \eqref{physical conc lp} we also have
$$\int_{|x-x(t)| \leq C(\eta_1)/N(t)}|u(t)|^3 dx \gtrsim c(\eta_1) N(t)^{-1}.$$
Combining this with \eqref{int bound high} and using the triangle inequality we find
\begin{equation*}
c(\eta_1) N(t)^{-1} \lesssim \int_{|x-x(t)| \leq
C(\eta_1)/N(t)}|P_{\geq N_*}u(t,x)|^3 dx.
\end{equation*}
Integrating this over $I_0$ and comparing with \eqref{int bound low} proves \eqref{integral bound on N}.
\end{proof}

We can now (finally!) conclude the contradiction argument.
It remains to prove $\|u\|_{L^6_{t,x}(I_0 \times \R^4)} \lesssim C(\eta_0,\eta_1,\eta_2,\eta_3)$, which we expect
since the bound
\eqref{length I bound} shows that the interval $I$ is not long enough to allow the $L^6_{t,x}$-norm
of $u$ to grow too large.

\begin{proposition}\label{lemma noncon imp spacetime bdd}
We have
$$\|u\|_{\ls (I_0 \times \R^4)} \lesssim C(\eta_0,\eta_1,\eta_2,\eta_3).$$
\end{proposition}

\begin{proof}
We normalize $N_{min} = 1$.  Let $\delta =
\delta(\eta_0,N_{max}) > 0$ be a small number to be chosen
momentarily.  Partition $I_0$ into $O(|I_0|/\delta)$ subintervals
$I_1,\dots,I_J$ with $|I_j| \leq \delta$.  Let $t_j \in I_j$.
Since $N(t_j) \leq N_{max}$ we have from Corollary \ref{lemma freq
loc}
\begin{equation}\label{high freq eta0 bound}
\|P_{\geq C(\eta_0)N_{max}}u(t_j)\|_{\ho_x} \leq \eta_0.
\end{equation}
Let $\tilde u(t) = e^{i(t-t_j)\Delta}P_{<C(\eta_0)N_{max}}u(t_j)$
be the free evolution of the low and medium frequencies of
$u(t_j)$. The above bound becomes $$\|u(t_j) - \tilde
u(t_j)\|_{\ho_x} \leq \eta_0.$$

But from Bernstein, Sobolev embedding, and \eqref{kinetic energy bound} we get
\begin{align*}
\|\tilde u(t) \|_{L^6_{x}} &\lesssim C(\eta_0,N_{max}) \|\tilde
u(t_j)\|_{L^4_{x}}\\
&\lesssim C(\eta_0, N_{max})\|\tilde u(t_j)\|_{\ho_x}\\
&\lesssim C(\eta_0,N_{max})
\end{align*}
for all $t \in I_j$,
so $$\|\tilde u \|_{L^6_{t,x}(I_j \times \R^4)} \lesssim C(\eta_0,
N_{max}) \delta^{1/6}.$$  Similarly we have
\begin{align*}
\|\nabla (|\tilde u(t)|^2 \tilde u(t))\|_{L^{4/3}_x} & \lesssim
\|\nabla \tilde u(t) \|_{L^4_x} \|\tilde u(t)\|_{L^4_x}^2\\
& \lesssim C(\eta_0, N_{max}) \|\tilde u(t) \|_{\ho_x}^3\\
& \lesssim C(\eta_0, N_{max})
\end{align*}
which shows $$\|\nabla(|\tilde u|^2 \tilde u)\|_{L^2_t L^{4/3}_x (I_j \times \R^4)}
\lesssim C(\eta_0, N_{max}) \delta ^{1/2}.$$ From these two
estimates, \eqref{kinetic energy bound}, and Lemma \ref{lemma
short time} with $e = -|\tilde u|^2 \tilde u$ we see that
$$\|u\|_{L^6_{t,x}(I_j \times \R^4)} \lesssim 1$$  if
$\delta$ is chosen small enough.  Summing these bounds in $j$ and
using \eqref{length I bound} we get $$\|u\|_{L^6_{t,x}(I_0 \times \R^4)}
\lesssim C(\eta_0, N_{max}) |I_0| \lesssim C(\eta_0,\eta_1,\eta_2,\eta_3).$$
\end{proof}

%%%%%%%%%%%%%%%%%%%%%%%%%%%%%%%%%%%%%%%%%%%%%%%%%%%%%%%%%%%%%%%%%%%%%%%%%%%%%%%%%%%%%%%%%%%
%
%
%                                   Section
%
%
%%%%%%%%%%%%%%%%%%%%%%%%%%%%%%%%%%%%%%%%%%%%%%%%%%%%%%%%%%%%%%%%%%%%%%%%%%%%%%%%%%%%%%%%%%%

\section{Remarks} We first note that as a consequence of the bound \eqref{l6 bounds}, one gets scattering, asymptotic completeness, and
uniform regularity. Indeed, we have

\begin{proposition}
Let $u_0$ have finite energy. Then there exist finite energy solutions $u_{\pm}(t,x)$ to the free Schr\"odinger equation
$(i\partial_t+\Delta)u_{\pm}=0$ such that
$$
\|u_{\pm}(t)-u(t)\|_{\dot{H}^1_x}\rightarrow 0
$$
as $t\rightarrow \pm\infty$. Furthermore, the maps $u_0\mapsto u_{\pm}(0)$ are homeomorphisms from $\dot{H}^1_x(\R^4)$
to itself. Finally, if $u_0 \in H^s_x$ for some $s\geq 1$, then $u(t)\in H^s_x$ for all times $t$, and one has the
uniform bounds
$$
\sup_{t\in \R}\|u(t)\|_{H^s_x}\leq C(E(u_0),s)\|u_0\|_{H^s_x}.
$$
\end{proposition}

\begin{proof}
We will only prove the statement for $u_{+}$, since the proof for $u_{-}$ follows similarly.
Let us first construct the scattering state $u_{+}(0)$. For $t>0$ define $v(t) = e^{-it\Delta}u(t)$. We will show
that $v(t)$ converges in $\dot{H}^1_x$ as $t\rightarrow \infty$, and define $u_{+}(0)$ to be that limit.
Indeed, from Duhamel's formula \eqref{duhamel} we have
\begin{align}\label{v}
v(t) = u_0 - i\int_{0}^{t} e^{-is\Delta}(|u|^2u)(s)ds.
\end{align}
Therefore, for $0<\tau<t$,
$$
v(t)-v(\tau)=-i\int_{\tau}^{t}e^{-is\Delta}(|u|^2u)(s)ds.
$$
Lemma \ref{lemma linear strichartz} and Lemma \ref{lemma trilinear strichartz} yield
\begin{align*}
\|v(t)-v(\tau)\|_{\dot{H}^1_x}
     &=\|e^{it\Delta}[v(t)-v(\tau)]\|_{\dot{H}^1_x} \\
     &\lesssim \|\nabla (|u|^2u)\|_{L^1_tL_x^2([\tau,t]\times\R^4)} \\
     &\lesssim \|u\|_{\dot{S}^1([\tau,t]\times\R^4)}^2 \|u\|_{L^6_{t,x}([\tau,t]\times\R^4)}.
\end{align*}
However, Lemma \ref{lemma persistence regularity} implies that $\|u\|_{\dot{S}^1}$ is finite, while the bound \eqref{l6 bounds}
 implies that for any $\eps>0$ there exists $t_{\eps}\in \R_{+}$ such that
$\|u\|_{L^6_{t,x}([t,\infty)\times\R^4)}\leq \eps$ whenever $t>t_{\eps}$. Hence,
\begin{center}
$\|v(t)-v(\tau)\|_{\dot{H}^1_x}\rightarrow 0 \quad$ as $t,\tau\rightarrow \infty$.
\end{center}
In particular, this implies that $u_{+}(0)$ is well defined. Also, inspecting \eqref{v} one easily sees that
\begin{align}
u_{+}(0)=u_0- i\int_{0}^{\infty}e^{-is\Delta}(|u|^2u)(s)ds
\end{align}
and thus
\begin{align}\label{u+}
u_{+}(t)=e^{it\Delta}u_0- i\int_{0}^{\infty}e^{i(t-s)\Delta}(|u|^2u)(s)ds.
\end{align}
By the same arguments as above, \eqref{u+} and Duhamel's formula \eqref{duhamel} imply that
$\|u_{+}(t)-u(t)\|_{\dot{H}^1_x}\rightarrow 0$ as $t\rightarrow\infty$.

Similar estimates prove that the map $u_0\mapsto u_{+}(0)$ is continuous from $\dot{H}^1_x(\R^4)$
to itself. To show that the map is also injective, let $u_{01}, u_{02}\in \dot{H}^1_x$, and let $u_1, u_2$
be the corresponding solutions to \eqref{schrodinger equation} with initial data $u_{01}$, $u_{02}$.
Assume that $\|e^{-it\Delta}u_j(t)-u_{+}(0)\|_{\dot{H}^1_x}\rightarrow 0$ as $t\rightarrow\infty$ for $j=1,2$. Then,
\eqref{u+} yields
\begin{align}
u_j(t)= u_{+}(t)- i\int_{t}^{\infty}e^{i(t-s)\Delta}(|u_j|^2u_j)(s)ds
\end{align}
for $t>0$ and $j=1,2$.  Using the Strichartz estimates of Lemma \ref{lemma linear strichartz}
and Lemma \ref{lemma trilinear strichartz}, Picard's fixed point argument forces $u_1(t)=u_2(t)$ for $t$ sufficiently
large. By the uniqueness of
solutions to the Cauchy problem \eqref{schrodinger equation} we obtain $u_{01}=u_{02}$, which proves injectivity for
the map $u_0\mapsto u_{+}(0)$ on $\dot{H}^1_x(\R^4)$.

We now turn towards the construction of the wave operators, i.e. for every $u_{+}(0)\in \dot{H}^1_x(\R^4)$ there exists
$u_0\in \dot{H}^1_x(\R^4)$ such that $\|u_{+}(t)-u(t)\|_{\dot{H}^1_x}\rightarrow 0$ as $t\rightarrow\infty$,
and moreover the map $u_{+}(0)\mapsto u_0$ from $\dot{H}^1_x(\R^4)$ to itself is continuous. For $t>0$
take $u_{+}(t)$ to be initial data for the equation \eqref{schrodinger equation} and solve the Cauchy problem backwards
in time. Denote the solution at time $t=0$ by $u_{0t}$. Choosing $t_2>t_1$ sufficiently large, we see that the free and
the nonlinear evolutions from $t_2$ to $t_1$ are within $\eps$ of each other in the $\dot{S}^1([t_1,t_2]\times\R^4)$-norm.
An application of Lemma \ref{lemma long time} with $e=0$ yields $\|u_{0t_1}-u_{0t_2}\|_{\dot{H}^1_x}\lesssim \eps$,
which implies that $u_{0t}$ converges in $\dot{H}^1_x(\R^4)$ as $t\rightarrow\infty$. Denote the limit by $u_0$, and let
$u$ be the global solution to the Cauchy problem \eqref{schrodinger equation} with initial data $u_0$. Then at time
$t>0$ we have
\begin{align*}
u(t) &= e^{it\Delta}u_0 - i\int_0^t e^{i(t-s)\Delta}(|u|^2u)(s)ds \\
     &= e^{it\Delta}u_{0t} + e^{it\Delta}(u_0-u_{0t}) - i\int_0^t e^{i(t-s)\Delta}(|u|^2u)(s)ds \\
     &= u_{+}(t) + e^{it\Delta}(u_0-u_{0t}),
\end{align*}
which implies that $\|u_{+}(t)-u(t)\|_{\dot{H}^1_x}\rightarrow 0$ as $t\rightarrow\infty$.

The continuity of the map $u_{+}(0)\mapsto u_0$ on $\dot{H}^1_x(\R^4)$ follows immediately from Lemma
\ref{lemma long time} with $e=0$.

The regularity statement follows from obvious modifications of Lemma \ref{lemma persistence regularity}.
\end{proof}

We now develop the tower-type bounds for $M(E)$.  To do so we must
determine the dependence of each of the $\eta_j$'s on the previous
ones.  Throughout, we will let $c$ and $C$ represent small and
large constants, possibly depending on the energy.  To begin, we need $\eta_0$ to be
small relative to the energy.  Specifically, we choose $\eta_0
\leq c E^{-C}$, so that when we remove the low- and high-frequency
portions of the solution (whose norms are bounded by $c \eta_0^C$) the majority of the energy remains.

Recall that the dependence of Lemma \ref{lemma long time} is
exponential in its parameters.  That is, if $E$, $E'$, and $M$ represent
the various bounds in the statement of the lemma, we need to
choose $\eps_1 \lesssim \exp (-M^C \langle E \rangle ^C \langle E'
\rangle ^C)$. When we apply this lemma in
the proof of Proposition \ref{lemma freq loc imp spacetime bds}, the parameters are
\begin{align*}
E'= &\eps \\
M= &\eta^{-C} M(\ecrit - \eta^C)^C \\
E= &\eta^{-C} M(\ecrit - \eta^C)^C \\
\eps_1> &\eta^{-C} M(\ecrit - \eta^C)^C \eps^{1/2},
\end{align*}
where the $\eta^{-C} M(\ecrit - \eta^C)^C$ factors arise from an application
of Lemma \ref{lemma induct on energy} earlier in the proof. Thus, we need to choose
$$\eta^{-C} M(\ecrit - \eta^C)^C \eps^{1/2} \lesssim \exp (-\eta^{-C} M(\ecrit -\eta^C)^C
\langle \eps \rangle^C),$$ or equivalently $$\eps^{\eps^{-C}} \lesssim \exp
(-C \eta^{-C} M(\ecrit - \eta^C)^C).$$ Also, in order
to make the pigeonhole argument work (still in the proof of
Proposition \ref{lemma freq loc imp spacetime bds}) we need the
frequency separation to be $$K(\eta) \geq \eps^{-\eps^{-2}} \geq C
\exp (C \eta^{-C}  M(\ecrit - \eta^C)^C).$$ As a
result, the dependence of the constants $C(\eta)$ and $c(\eta)$ in
Corollary \ref{lemma freq loc} is quite bad:
\begin{gather*}
C(\eta) \geq C \exp (C \eta^{-C} M(\ecrit - \eta^C)^C)\\
c(\eta) \leq 1/C(\eta).
\end{gather*}
This motivates introducing the notation $t(\eta)$ for any term
of the form $$c \exp \bigl(-C \eta^{-C} M(\ecrit - \eta^C)^C\bigr).$$

Now, the constant $c(\eta_0)$ appearing in the proof of
Proposition \ref{lemma potential bdd below} needs to be of the
form $c(\eta_0) = t(\eta_0)$, due to the use of Corollary
\ref{lemma freq loc}. So in order to apply Lemma \ref{lemma long time}
(later in the same proof) we need to choose $\eta_1 \leq t(c(\eta_0)) =
t(t(\eta_0))$.

The constant $c(\eta_1)$ appearing in Proposition \ref{prop flim} is roughly of size $t(\eta_1)$.
In order to apply this proposition in Section 6, we therefore must choose $\eta_2 \leq t(\eta_1)$.

We need to choose $\eta_3$ smaller than any polynomial in $\eta_2$ appearing
in the proof of Proposition \ref{energyevac}.  For instance, $\eta_3 \leq \exp (- \eta_2^{-C})$ suffices;
given the final bound we obtain, this makes $\eta_2$ and $\eta_3$ virtually equivalent.

We also see from the proof of Proposition \ref{energyevac} that $N_{max} \leq
1/t(\eta_3)$. Proposition \ref{lemma noncon imp spacetime bdd} then implies the desired contradiction,
provided we choose $\eta_4 \leq t(\eta_3)$.

Putting this all together we find a final bound of the form
\begin{equation}\label{first tower}
M(E) \leq 1/t(t(t(E^{-C}))).
\end{equation}

To properly simplify this expression requires some notation (the Ackerman
hierarchy).  Recall that multiplication is iterated addition
$$a \times b = a + \dots + a$$
with $b$ factors on the right, and exponentiation is iterated
multiplication
$$a \uparrow b = a \times \dots \times a.$$
We define tower exponentiation
as iterated exponentiation:
$$a \uparrow \uparrow b := a \uparrow (a \uparrow \dots (a
\uparrow a) \dots )$$ with $b$ arrows on the right.  Similarly,
double tower exponentiation is iterated tower exponentiation
$$a \uparrow \uparrow \uparrow b := a \uparrow \uparrow (a \uparrow
\uparrow \dots (a
\uparrow \uparrow a) \dots ),$$ and triple tower exponentiation is
iterated double tower exponentiation, etc.

Now, we can essentially expand \eqref{first tower} as
\begin{equation}\label{second tower}
M(E) \leq \exp (M(E - t(t(E^{-C})))).
\end{equation}
Iterating \eqref{second tower} approximately $1/t(t(E^{-C}))$ times, we
obtain
$$M(E) \leq C \uparrow \uparrow 1/t(t(E^{-C})).$$
We again expand this as
$$M(E) \leq C \uparrow \uparrow \exp (M(E - t(E^{-C}))),$$
and iterating it we obtain
$$M(E) \leq C \uparrow \uparrow \uparrow 1/t(E^{-C}).$$
Repeating this process one more time we get the final bound
$$M(E) \leq C \uparrow \uparrow \uparrow \uparrow CE^C.$$

This tower bound is a far cry from the exponential bound
obtained in \cite{tao: gwp radial}, though better than that obtained in \cite{ckstt:gwp}, which is
of the form $M(E) \leq C \uparrow \uparrow \uparrow \uparrow
\uparrow \uparrow \uparrow \uparrow CE^C$.  It is plausible that
one could obtain a polynomial bound in the
energy-critical case, but to do so would require abandoning
the inductive approach used in this paper.

\end{document}